\documentclass{article}
\usepackage{amsfonts}

\usepackage{amsmath}
\usepackage{scalefnt}

\newtheorem{theorem}{Theorem}[section]

\newtheorem{corollary}[theorem]{Corollary}

\newtheorem{lemma}[theorem]{Lemma}

\newtheorem{proposition}[theorem]{Proposition}

\newenvironment{proof}[1][Proof]{\noindent\textbf{#1.} }{\ \rule{0.5em}{0.5em}}

\input{tcilatex}

\begin{document}

\title{On the Blaschke Conjecture for $3$-Webs}
\author{ Vladislav V. Goldberg and Valentin V. Lychagin}
\maketitle

\begin{abstract}
We find relative differential invariants of orders eight and nine for a
planar nonparallelizable $3$-web such that their vanishing is necessary and
sufficient for a $3$-web to be linearizable. This solves the Blaschke
conjecture for $3$-webs. As a side result, we show that the number of
linearizations in the Gronwall conjecture does not exceed fifteen and give
criteria for rigidity of $3$-webs.


\end{abstract}

\section{Introduction}

Let $W_{d}$ be a $d$-web given by $d$ one-parameter foliations of curves on
a two-dimensional manifold $M^{2}$. The web $W_{d}$ is linearizable
(rectifiable) if it is equivalent to a linear $d$-web, i.e., a $d$-web
formed by $d$ one-parameter foliations of straight lines on a projective
plane.

The problem of finding a criterion of linearizability \ of webs was posed by
Blaschke in the 1920s (see, for example, his book \cite{B 55}, $\S 17$ and $%
\S 42$) who claimed that it is hopeless to find such a criterion$.$\
Comparing\ the numbers of relative invariants for a general $3$-web $W_{3}$
(and a general $4$-web $W_{4}$) and a linear $3$-web (and a linear $4$-web),
Blaschke made the conjectures that conditions of linearizability for a $3$%
-web $W_{3}$ should consist of four relations for the ninth order web
invariants (four PDEs of ninth order) and those for a $4$-web $W_{4}$ should
consist of two relations for the fourth order web invariants (two PDEs of
fourth order) .

In \cite{AGL} the authors proved that the Blaschke conjecture on
linearizability conditions for $4$-webs was correct: a $4$-web $W_{4}$ is
linearizable if and only if its two fourth order invariants vanish.\ In \cite%
{AGL} a complete solution of the linearizability problem for $d$-webs, $%
d\geq 5,$ was also presented. In \cite{G 04} the linearizability conditions
found in \cite{AGL} were applied to check whether some known classes of $4$%
-webs are linearizable.

In the present paper we\ continue to use the Akivis approach (see \cite{AGL}%
) for establishing criteria of linearizability of $3$-webs. In this
approach, the linearizability problem is reduced to the solvability of the
system of nonlinear partial differential equations on the components of the
affine deformation tensor. This is the system of four nonlinear first-order
PDEs on three functions defined on the plane. In the paper \cite{G 89} the
first obstruction for integrability of the system was found. In this paper
we use results of \cite{KL} to investigate the integrability of the system
and show that the obstruction found in \cite{G 89} coincides with the Mayer
bracket defined in \cite{KL}.

We show that for nonparallelizable $3$-webs, the solvability of the system
indicated above is equivalent to the existence of real and smooth solutions
of the system of five algebraic equations of degrees not exceeding $17,18,18$
and $24,24.$ This allows us:

\begin{enumerate}
\item[(i)] To find relative differential invariants whose vanishing leads to
the linearizability of a $3$-web $W_{3}$. This solves the \emph{Blaschke
problem} mentioned earlier on finding linearizability conditions in the form
of invariants whose vanishing is necessary and sufficient for
linearizability of a $3$-web $W_{3}$. There are two types of invariants: $18$
of them have order eight and $1040$ have order nine. Note that the number of
invariants can be different but there are always invariants of order eight.
Note also that the Blaschke estimation of the ''functional codimension''\ of
the orbits of the linearizable $3$-webs was correct, but the number of
invariants was not. Moreover, the problem has invariants of order eight that
do not match his prediction.

\item[(ii)] To establish the algorithm for determining whether a given $3$%
-web $W_{3}\;$is linearizable. This algorithm is based on investigation of
the existence of a real solution of the five algebraic equations mentioned
above.
\end{enumerate}

We have checked that the differential invariants vanish for all linear $3$%
-webs $W_{3}\;$and apply the algorithm to two more examples (of nonlinear) $%
3 $-webs $W_{3}.$

As a side result, we obtain an estimation for the Gronwall conjecture. \ In
1912 Gronwall (\cite{gr})\ made the following conjecture: \textit{if a
nonparallelizable }$3$\textit{-web }$W_{3}$\textit{\ in the plane is
linearizable, then, up to a projective transformation, a diffeomorphism
transforming }$W_{3}$\textit{\ into a linear }$3$\textit{-web is uniquely
determined. ~}The Gronwall conjecture is also called the ''fundamental
theorem''\ of nomography. Note that for parallelizable $3$-webs such
uniqueness does not take place. In fact, such a $3$-web \ is formed by the
tangents to a curve of third degree, but curves of third degree have
nontrivial projective invariants (see \cite{B 55}, \S 17).

Bol (\cite{bol1}, \cite{bol2}, 1938) and Bor\r{u}vka (\cite{bo}, 1938)
proved that the number of projectively nonequivalent linearizations of a
nonparallelizable, linearizable $3$-web does not exceed $16.\ $Grifone,\
Muzsnay and Saab (\cite{GMS 01}, 2001) proved that this number does not
exceed $15.$ We also prove that this number does not exceed $15,$ and give
criteria for rigidity of $3$-webs, but our method is different from that in %
\cite{GMS 01}.

Note that Vaona (\cite{v2}, 1961) and Smirnov (\cite{sm1}, \cite{sm2})
considered the Gronwall conjecture from the point of view of nomography.
Vaona claimed that the\ above mentioned number does not exceed $11,\;$and
Smirnov claimed that this number does not exceed one (i.e., that the
Gronwall conjecture is right).\ 

In addition, we find the linearity condition for $3$-webs and establish the
relationship of this to the condition that a plane curve consists of flexes$%
\;$and to the Euler equation in gas-dynamics.

The completion of this paper would not have been possible without the
support provided to the authors by the Mathematisches Forschungsinstitut
Oberwolfach (MFO), Germany. We express our deep gratitude to Professor Dr.
G.-M. Greuel, the director of MFO, for the opportunity to use the excellent
facilities at MFO.

\section{Basics Constructions}

We recall main constructions for $3$-webs on two-dimensional manifolds (see,
for example, \cite{BB 38} or \cite{B 55}, \cite{G 89}) in a form suitable
for us.

Let $M^{2}$ be a two-dimensional manifold, and suppose that a $3$-web $W_{3}$
is given on $M^{2}$ by three differential $1$-forms $\omega _{1},\omega
_{2}, $ and $\omega _{3}$ such that any two of them are linearly independent.

\begin{proposition}
The forms $\omega _{1},\omega _{2},$ and $\omega _{3}$ can be normalized in
such a way that the normalization condition 
\begin{equation}
\omega _{1}+\omega _{2}+\omega _{3}=0  \label{normalization equation}
\end{equation}%
holds.
\end{proposition}

\begin{proof}
In fact, if we take the forms $\omega _{1}$ and $\omega _{2}$ as co-basis
forms of $M^{2}$, then the form $\omega _{3}$ is a linear combination of the
forms $\omega _{1}$ and $\omega _{2}$: 
\begin{equation*}
\omega _{3}=\alpha \omega _{1}+\beta \omega _{2\,},
\end{equation*}%
where $\alpha ,\beta \neq 0.\;$After the substitution 
\begin{equation*}
\omega _{1}\rightarrow \frac{1}{\alpha }\omega _{1},\;\omega _{2}\rightarrow 
\frac{1}{\beta }\omega _{2},\;\omega _{3}\rightarrow -\omega _{3}
\end{equation*}%
the above equation becomes (\ref{normalization equation}).
\end{proof}

It is easy to see that any two of such normalized triplets $\omega
_{1},\omega _{2},\omega _{3}$ and $\omega _{1}^{s},\omega _{2}^{s},\omega
_{3}^{s}$ determine the same $3$-web $W_{3}$ if and only if 
\begin{equation}
\omega _{1}^{s}=s^{-1}\omega _{1},\ \omega _{2}^{s}=s^{-1}\omega _{2},\
\omega _{3}^{s}=s^{-1}\omega _{3}  \label{perenormirovka}
\end{equation}%
for a non-zero smooth function $s\in C^{\infty }\left( M^{2}\right) .$

\subsection{Structure Equations}

From now on we shall assume that a $3$-web $W_{3}$ is given by differential $%
1$-forms $\omega _{1},\omega _{2},$ and $\omega _{3}$ normalized by
condition (\ref{normalization equation}).

Because $M^{2}$ is\ a two-dimensional manifold, there is a unique
differential $1$-form $\gamma $ such that 
\begin{equation}  \label{web structure equations}
\renewcommand{\arraystretch}{1.3} 
\begin{array}{ll}
d\omega _{1}=\omega _{1}\wedge \gamma , &  \\ 
d\omega _{2}=\omega _{2}\wedge \gamma . & 
\end{array}
\renewcommand{\arraystretch}{1}
\end{equation}

Moreover, it follows \ from (\ref{normalization equation}) that 
\begin{equation*}
d\omega _{3}=\omega _{3}\wedge \gamma .
\end{equation*}

We call $\gamma $ the \emph{connection form }and equations (\ref{web
structure equations}) the \emph{web structure equations.}

Later on we shall see that $\gamma $ determines the so-called Chern
connection on $M^{2}.$

For other representations $\left( \omega _{1}^{s},\omega _{2}^{s},\omega
_{3}^{s}\right) $ of the web, structure equations (\ref{web structure
equations}) take the form

\begin{eqnarray*}
d\omega _{1}^{s} &=&\omega _{1}^{s}\wedge \gamma ^{s}, \\
d\omega _{2}^{s} &=&\omega _{2}^{s}\wedge \gamma ^{s},
\end{eqnarray*}%
where%
\begin{equation*}
\gamma ^{s}=\gamma +\frac{ds}{s}.
\end{equation*}%
Note that the differential $2$-form $d\gamma $ does not depend on the web
representation and is an invariant of $3$-webs.

Let 
\begin{equation*}
d\gamma ^{s}=K_{s}~\omega _{1}^{s}\wedge \omega _{2}^{s}
\end{equation*}%
and 
\begin{equation*}
d\gamma =K~\omega _{1}\wedge \omega _{2}.
\end{equation*}%
The function $K$ is called the \emph{web curvature. }It follows from the
last two equations that 
\begin{equation*}
K_{s}=s^{2}K.
\end{equation*}%
This means that the web curvature $K$ is a relative invariant of weight two.

Let $\partial _{1},\partial _{2}$ be the dual basis of the vector field
module: $\omega _{i}\left( \partial _{j}\right) =\delta _{ij},$ $i,j=1,2.$
One has 
\begin{equation*}
df=\partial _{1}\left( f\right) ~\omega _{1}+\partial _{2}\left( f\right)
~\omega _{2}
\end{equation*}%
for smooth functions $f\in C^{\infty }\left( M^{2}\right) .$

If we decompose the connection forms $\gamma $ and $\gamma ^{s}$ relative to
the basis $\{\omega _{1},\omega _{2}\}$: 
\begin{equation}
\gamma =g_{1}~\omega _{1}+g_{2}~\omega _{2}  \label{gamma}
\end{equation}%
and%
\begin{equation*}
\gamma ^{s}=g_{s1}~\omega _{1}^{s}+g_{s2}~\omega _{2}^{s},
\end{equation*}%
we get%
\begin{eqnarray*}
g_{s1} &=&sg_{1}+\partial _{1}s, \\
g_{s2} &=&sg_{2}+\partial _{2}s.
\end{eqnarray*}%
In addition, we find%
\begin{equation}
\lbrack \partial _{1},\partial _{2}]=-g_{2}\partial _{1}+g_{1}\partial _{2}.
\label{commutator}
\end{equation}%
This follows from%
\begin{equation*}
\omega _{1}\left( [\partial _{1},\partial _{2}]\right) =-d\omega _{1}\left(
\partial _{1},\partial _{2}\right) =(\gamma \wedge \omega _{1})\left(
\partial _{1},\partial _{2}\right) =-\gamma \left( \partial _{2}\right)
=-g_{2}
\end{equation*}%
and 
\begin{equation*}
\omega _{2}\left( [\partial _{1},\partial _{2}]\right) =-d\omega _{2}\left(
\partial _{1},\partial _{2}\right) =(\gamma \wedge \omega _{2})\left(
\partial _{1},\partial _{2}\right) =\gamma \left( \partial _{1}\right)
=g_{1}.
\end{equation*}%
Remark that 
\begin{equation*}
\gamma \left( \lbrack \partial _{1},\partial _{2}]\right) =0.
\end{equation*}%
For the curvature function, one has%
\begin{equation}
K=\partial _{1}\left( g_{2}\right) -\partial _{2}\left( g_{1}\right) ,
\label{curvature equation for web}
\end{equation}%
because%
\begin{eqnarray*}
d\gamma &=&dg_{1}\wedge \omega _{1}+dg_{2}\wedge \omega _{2}+g_{1}d\omega
_{1}+g_{2}d\omega _{2}= \\
&&-\partial _{2}\left( g_{1}\right) ~\omega _{1}\wedge \omega _{2}+-\partial
_{1}\left( g_{2}\right) ~\omega _{1}\wedge \omega _{2}+g_{1}\omega
_{1}\wedge \gamma +g_{2}\omega _{2}\wedge \gamma \\
&=&-\partial _{2}\left( g_{1}\right) ~\omega _{1}\wedge \omega
_{2}+-\partial _{1}\left( g_{2}\right) ~\omega _{1}\wedge \omega
_{2}+g_{1}g_{2}\omega _{1}\wedge \omega _{2}-g_{1}g_{2}\omega _{1}\wedge
\omega _{2} \\
&=&\left( \partial _{1}\left( g_{2}\right) -\partial _{2}\left( g_{1}\right)
\right) \omega _{1}\wedge \omega _{2}.
\end{eqnarray*}

In this paper we shall apply the following two normalizations: (i) $d\omega
_{3}=0,$ and (ii) $K=1.$

The first one defines a $3$-web up to gauge transformations: $f\rightarrow
F\left( f\right) ,$ while the second one defines the $e$-structure on $%
M^{2}. $

Below we consider these two normalizations in detail.

\subsection{Normalization $\mathbf{\emph{d}\protect\omega }_{3}=\mathbf{0}$}

We assume that $M^{2}$ is a simply connected domain of $\mathbb{R}^{2}$, and
therefore there exists a smooth function $f$ such that $\omega _{3}$ is
proportional to $df,$ that is, $\omega _{3}\wedge df=0.$ The function $f$ is
called the \emph{web function. }

Note that this function is defined up to a renormalization (gauge
transformation) $f\longmapsto F\left( f\right) .$

We choose a representation of $W_{3}$ such that 
\begin{equation}
\omega _{3}=df.  \label{normalized condition2}
\end{equation}

Similarly, one finds smooth functions $x$ and $y$ for forms $\omega _{1}$
and $\omega _{2}$ such that%
\begin{equation*}
\omega _{1}=a\;dx,\ \omega _{2}=b\;dy
\end{equation*}%
for some smooth functions $a$ and $b.$

Moreover, the functions $x$ and $y$ are independent and therefore can be
viewed as (local) coordinates. In these coordinates, the normalization
condition gives 
\begin{equation*}
\omega _{1}=-f_{x}\;dx,\ \omega _{2}=-f_{y}\;dy,\ \omega _{3}=df.
\end{equation*}%
The vector fields $\partial _{1}$ and $\partial _{2}$ \ take the following
form 
\begin{equation*}
\partial _{1}=-\frac{1}{f_{x}}\frac{\partial }{\partial x},\ \ \ \partial
_{2}=-\frac{1}{f_{y}}\frac{\partial }{\partial y}.
\end{equation*}

In this case%
\begin{equation*}
0=d\omega _{3}=\omega _{3}\wedge \gamma
\end{equation*}%
and 
\begin{equation*}
\gamma =-H\omega _{3}=H\left( \omega _{1}+\omega _{2}\right)
\end{equation*}%
for some function $H.$

Hence (see (\ref{gamma})) 
\begin{equation*}
g_{1}=g_{2}=H.
\end{equation*}

In terms of the web function $f$, one has 
\begin{equation*}
H=\frac{f_{xy}}{f_{x}f_{y}},
\end{equation*}%
and 
\begin{equation*}
\gamma =-\frac{f_{xy}}{f_{x}f_{y}}\omega _{3}.
\end{equation*}

For the curvature function $K$ one gets the following expression: 
\begin{equation*}
K=-\frac{1}{f_{x}f_{y}}\left( \log \left( \frac{f_{x}}{f_{y}}\right) \right)
_{xy}=\frac{f_{xyy}}{f_{x}f_{y}^{2}}-\frac{f_{xxy}}{f_{x}^{2}f_{y}}+\frac{%
f_{xx}f_{xy}}{f_{x}^{3}f_{y}}-\frac{f_{xy}f_{yy}}{f_{x}f_{y}^{3}}
\end{equation*}%
(cf. \cite{B 55}, \ \S\ 9, or \cite{AS 92}, p. 43).

For the basis vector fields $\partial _{1}$ and $\partial _{2}$, the
structure equations take the form 
\begin{equation}
\lbrack \partial _{1},\partial _{2}]=H~(\partial _{2}-\partial _{1}),
\label{vector structure eq}
\end{equation}%
and 
\begin{equation}
K=\partial _{1}\left( H\right) -\partial _{2}\left( H\right) .
\label{curvature main}
\end{equation}

\subsection{Normalization $\emph{K}\mathbf{=1}$}

In this section we assume that $K$ is a nonvanishing function: $K\neq 0.$ We
can assume that $K>0$ (changing the orientation if necessary), that is, 
\begin{equation*}
K=k^{2}
\end{equation*}%
for some weight one smooth function $k.$

Let us take $s=k^{-1}$ and denote by $\theta _{i}$ the differential $1$%
-forms $\omega _{i}^{s}$ with $s=k^{-1}$:%
\begin{equation*}
\theta _{i}=k\omega _{i}
\end{equation*}%
for $i=1,2.$

We shall denote the corresponding connection form $\gamma ^{s}$ by $\alpha $:%
\begin{equation*}
\alpha =\gamma -\frac{dk}{k}.
\end{equation*}%
One has $k_{t}=tk$ for any positive smooth function $t,\ $and therefore$\
\theta _{i}=k\omega _{i}=k_{t}\omega _{i}^{t},$ $i=1,2,$ are invariant
differential $1$-forms intrinsically connected with the web. They define the 
$e$-structure on $M^{2}$ and satisfy the structure equations 
\begin{equation}
\renewcommand{\arraystretch}{1.3}%
\begin{array}{ll}
d\theta _{1}=\theta _{1}\wedge \alpha , &  \\ 
d\theta _{2}=\theta _{2}\wedge \alpha , &  \\ 
d\alpha =\theta _{1}\wedge \theta _{2}, & 
\end{array}%
\renewcommand{\arraystretch}{1}  \label{standard structure equations}
\end{equation}%
because $K_{k^{-1}}=(k^{-1})^{2}K=1.$

Let $\{\nabla _{1},\nabla _{2}\}$ be the basis dual to the co-basis $\left\{
\theta _{1},\theta _{2}\right\} ,$ and let 
\begin{equation*}
\alpha =a_{1}~\theta _{1}+a_{2}~\theta _{2}.
\end{equation*}%
Then (\ref{commutator}) and (\ref{curvature equation for web}) imply that 
\begin{equation}
\lbrack \nabla _{1},\nabla _{2}]=-a_{2}~\nabla _{1}+a_{1}~\nabla _{2}
\label{nabla commutator}
\end{equation}%
and%
\begin{equation}
\nabla _{1}\left( a_{2}\right) -\nabla _{2}\left( a_{1}\right) =1,
\label{normalized K}
\end{equation}%
where $a_{1}$ and $a_{2}$ are invariants of the web.

In terms of the web function $f,$ one has$\ $%
\begin{equation}
a_{1}=\frac{H}{k}-\frac{\partial _{1}k}{k^{2}},\ \ \ a_{2}=\frac{H}{k}-\frac{%
\partial _{2}k}{k^{2}}.  \label{a-invariants}
\end{equation}

\subsection{Linear 3-Webs}

In this section we consider linear $3$-webs. Let $W_{3}\;$be a $3$-web given
by a web function $z=f(x,y).$ The following theorem gives us a criterion for 
$W_{3}$ to be linear.

\begin{theorem}
Suppose that a $3$-web $W_{3}$\ is given locally by the function $z=f(x,y).$
Then $W_{3}$\textit{\ is linear if and only if}%
\begin{equation}
f_{y}^{2}~f_{xx}-2f_{x}f_{y}~f_{xy}+f_{x}^{2}~f_{yy}=0.  \label{lincondition}
\end{equation}
\end{theorem}

\begin{proof}
Note that a 3-web $W_{3}\;$can be also given by a nonvanishing function $%
f_{x}(x,y)/f_{y}\left( x,y\right) .\;$Namely, the horizontal and vertical
leaves are given by $x=\limfunc{const}$ and $y=\limfunc{const},$
respectively, and the transversal leaves are defined in such a way that $%
t=\tan \alpha ,\,$where $\alpha $ is the angle of the normal to the
transversal leaves with the horizontal leaves. So, the web $W_{3}\;$is
linear if and only if the function $\,f_{x}(x,y)/f_{y}\left( x,y\right) $\
remains constant along the transversal leaves. Thus%
\begin{equation*}
d\left( \frac{f_{x}}{f_{y}}\right) =0\;\func{mod}\;(\omega _{1}+\omega _{2})
\end{equation*}%
and%
\begin{equation*}
\partial _{1}\left( \frac{f_{x}}{f_{y}}\right) \omega _{1}+\partial
_{2}\left( \frac{f_{x}}{f_{y}}\right) \omega _{2}=0\;\func{mod}\;(\omega
_{1}+\omega _{2})
\end{equation*}%
or 
\begin{equation}
\partial _{1}\left( \frac{f_{x}}{f_{y}}\right) -\partial _{2}\left( \frac{%
f_{x}}{f_{y}}\right) =0.  \label{lincondition3}
\end{equation}%
It is easy to see that equation (\ref{lincondition3}) is equivalent to
equation (\ref{lincondition}).
\end{proof}

\textbf{Remark.} Note that linearity condition $(\ref{lincondition})$ of a $%
3 $-web $W_{3}\;$can be written in the determinant form:%
\begin{equation}
\ \ \det \left\Vert 
\begin{array}{ccc}
f_{xx} & f_{xy} & f_{x} \\ 
f_{xy} & f_{yy} & f_{y} \\ 
f_{x} & f_{y} & 0%
\end{array}%
\right\Vert \ =0.  \label{lincondition2}
\end{equation}

Note also that linearity condition $(\mathrm{\ref{lincondition}})$ $($or $(%
\mathrm{\ref{lincondition2}}))$ for a $3$-web is also the necessary and
sufficient condition for a point $(x,y)$ to be a flex of the curve defined
by the equation $f(x,y)=0\;($see, for example, $\cite{sh}$, section \textrm{%
1.1.5}$)$. The difference is that here $(\mathrm{\ref{lincondition}})$ is
the equation for finding the function $z=f(x,y)\;$(it should be satisfied
for all points $(x,y)$) while in algebraic geometry $(\mathrm{\ref%
{lincondition}})$ is the equation for finding the flexes $(x,y)\;$of the
curve defined by the equation $f(x,y)=0\;$provided that the function $%
f(x,y)\;$is given.

Differential equation (\ref{lincondition}) can be integrated as follows. Let
us rewrite this equation in form (\ref{lincondition3}). Then 
\begin{equation*}
\partial _{x}\left( \frac{f_{x}}{f_{y}}\right) -\left( \frac{f_{x}}{f_{y}}%
\right) \partial _{y}\left( \frac{f_{x}}{f_{y}}\right) =0,
\end{equation*}%
or setting 
\begin{equation*}
w=\frac{f_{x}}{f_{y}},
\end{equation*}%
we can rewrite (\ref{lincondition}) as the following system:%
\begin{eqnarray*}
\partial _{x}w-w\partial _{y}w &=&0, \\
\partial _{x}f-w\partial _{y}f &=&0.
\end{eqnarray*}

The first equation 
\begin{equation*}
\partial _{x}w-w\partial _{y}w=0
\end{equation*}%
is the Euler equation in gas-dynamics (see, for example, \cite{LL}, p. 3).

Solutions of this equation are well-known. Namely, if $w_{0}\left( y\right)
=\left. w\right| _{x=0}$ gives a Cauchy data, then the solution $w(x,y)$ can
be found from the system 
\begin{equation}  \label{EulerWeb}
\renewcommand{\arraystretch}{1.3} 
\begin{array}{lll}
y+w_{0}\left( \lambda \right) x-\lambda & = & 0, \\ 
w(x,y)-w_{0}\left( \lambda \right) & = & 0%
\end{array}
\renewcommand{\arraystretch}{1}
\end{equation}%
by elimination of the parameter $\lambda .$

Further, if $w$ is a solution of the Euler equation, then the functions $w$
and $f$ are first integrals of the vector field 
\begin{equation*}
\partial _{x}-w\partial _{y},
\end{equation*}%
and therefore there is the relation $f=F\left( w\right) $ for some smooth
function $F.$

Summarizing we get the following description of linear $3$-webs.

\begin{proposition}
The web functions $f\left( x,y\right) $ of linear $3$-webs have the form 
\begin{equation*}
f\left( x,y\right) =F\left( w\left( x,y\right) \right) ,
\end{equation*}%
where $w\left( x,y\right) $ is a solution of the Euler equation, and $F\;$is
some smooth function.
\end{proposition}

As we saw earlier, the web functions are defined up to gauge transformations 
$f\longmapsto F\left( f\right) .$ Therefore, the above proposition yields
the following description of linear $3$-webs.

\begin{theorem}
Web functions of linear $3$-webs can be chosen as solutions of the Euler
equation.
\end{theorem}

\textbf{Example 1} Taking $w_{0}\left( y\right) =y,$\ we get the linear $3$%
-web with the web function $w=y/\left( 1-x\right) .\;$This $3$-web is
generated by two families of coordinate lines $\left\{ x=\limfunc{const}%
\right\} $, $\left\{ y=\limfunc{const}\right\} $ and the straight lines of
the pencil with the center $(1,0).\;$This $3$-web is parallelizable.

\bigskip

\textbf{Example 2 }Taking $w_{0}\left( y\right) =y^{2}/4,$ \ we get the
linear $3$-web with the web function $\left( \frac{1+\sqrt{1-xy}}{x}\right)
^{2},$ or simply%
\begin{equation*}
f=\frac{1+\sqrt{1-xy}}{x}.
\end{equation*}%
It is easy to prove that this $3$-web is generated by two families of
coordinate lines $\left\{ x=\limfunc{const}\right\} $, $\left\{ y=\limfunc{%
const}\right\} $ and the tangents to the hyperbola $y=\frac{1}{x}.$ In fact,
the leaves of the third foliation of this web are level sets of the above
web function, i.e., they are determined by the equation%
\begin{equation*}
\frac{1+\sqrt{1-xy}}{x}=C,
\end{equation*}%
where $C\;$is a constant. The latter equation is equivalent to the equation%
\begin{equation*}
y=-C^{2}x+2C.
\end{equation*}%
Thus the leaves of the the third foliation are straight lines. To find the
envelope of these leaves, we differentiate the above equation with respect
to $C.\;$This gives $C=\frac{1}{x}.\;$Therefore, the envelope is defined by
the equation $y=\frac{1}{x}.$

\bigskip

\textbf{Example 3 }Taking $w_{0}\left( y\right) =-2\sqrt{-y},$ \ we get the
linear $3$-web with the web function%
\begin{equation*}
f=x+\sqrt{x^{2}-y}.
\end{equation*}%
Using the same approach as in Example 2, we can prove that the leaves of the
third foliation are straight lines defined by the equation%
\begin{equation*}
y=2Cx-C^{2},
\end{equation*}%
and these straight lines are tangent to the parabola $y=x^{2}.$

\section{The Chern connection}

Recall that a connection $\nabla $ in a vector bundle $\pi :E\left( \pi
\right) \rightarrow B$ over a manifold $B$ can be defined by a covariant
differential $d_{\nabla }:\Gamma \left( \pi \right) \rightarrow \Gamma
\left( \pi \right) \otimes \Omega ^{1}\left( B\right) ,$ where $\Gamma
\left( \pi \right) $ is the module of smooth sections of the bundle $\pi ,$
and $\Omega ^{1}\left( B\right) $ is the module of smooth differential $1$%
-forms on the manifold $B.$ The covariant differential can be extended in a
natural way to the following sequence: 
\begin{equation*}
\Gamma \left( \pi \right) \overset{d_{\nabla }}{\rightarrow }\Gamma \left(
\pi \right) \otimes \Omega ^{1}\left( B\right) \overset{d_{\nabla }}{%
\rightarrow }\Gamma \left( \pi \right) \otimes \Omega ^{2}\left( B\right) 
\overset{d_{\nabla }}{\rightarrow }\cdots
\end{equation*}%
The square of the covariant differential is the module homomorphism%
\begin{equation*}
d_{\nabla }^{2}\overset{\text{def}}{=}R_{\nabla }:\Gamma \left( \pi \right)
\rightarrow \Gamma \left( \pi \right) \otimes \Omega ^{2}\left( B\right) .
\end{equation*}%
This homomorphism $R_{\nabla }\ $is called the \emph{curvature }of the
connection $\nabla $.

We shall apply this construction to $3$-webs on a two-dimensional manifold $%
M.$ Let $\pi =\tau ^{\ast }:T^{\ast }\left( M\right) \rightarrow M$ be the
cotangent bundle, and let $W_{3}$ be a $3$-web defined by the differential $%
1 $-forms $\left\{ \omega _{1},\omega _{2},\omega _{3}\right\} $ normalized
by (\ref{normalization equation}).

We use the differential $1$-form $\gamma $ to define a connection in the
cotangent bundle by the following covariant differential: 
\begin{equation*}
d_{\gamma }:\Omega ^{1}\left( M\right) \rightarrow \Omega ^{1}\left(
M\right) \otimes \Omega ^{1}\left( M\right) ,
\end{equation*}%
where 
\begin{eqnarray*}
d_{\gamma }\left( \omega _{1}\right) &&\!\!\!\!=\;\;\!\!\!\!\;-\omega
_{1}\otimes \gamma , \\
d_{\gamma }\left( \omega _{2}\right) &&\!\!\!\!=\!\!\!\!\;\;\;-\omega
_{2}\otimes \gamma ;
\end{eqnarray*}%
and $\otimes $ denotes the tensor product.

Note that in the tensor product $\Omega ^{1}\left( M\right) \otimes \Omega
^{1}\left( M\right) $ the first factor plays the role of coefficients and
should be differentiated due to the connection, and the second one is
differentiated by the de Rham differential.

It is easy to check that the curvature form of the above connection is equal
to $-d\gamma ,$ that is, $d_{\gamma }^{2}:\Omega ^{1}\left( M\right)
\rightarrow \Omega ^{1}\left( M\right) \otimes \Omega ^{2}\left( M\right) $
is the multiplication by $-d\gamma $:%
\begin{equation*}
d_{\gamma }^{2}\left( \omega \right) =-\omega \otimes d\gamma
\end{equation*}%
for any differential form $\omega \in \Omega ^{1}\left( M\right) .$

This connection is called the \emph{Chern connection} of the web.

It is also easy to check that the Chern connection satisfies the relations 
\begin{equation*}
d_{\gamma ^{s}}\left( \omega _{i}^{s}\right) =-\omega _{i}^{s}\otimes \gamma
^{s}
\end{equation*}%
for $i=1,2,$ and any non-zero smooth function $s.$

The straightforward computation shows also that $d_{\gamma }$ is a
torsion-free connection.

Note that in the case $K\neq 0$ the second normalization $\left( K=1\right) $
leads us to the invariant $1$-forms $\theta _{1}$ and $\theta _{2}$ and to
the unique Chern connection $d_{\alpha }.$

Recall that for the covariant differential $d_{\nabla }:\Omega ^{1}\left(
M\right) \rightarrow \Omega ^{1}\left( M\right) \otimes \Omega ^{1}\left(
M\right) $ of any torsion-free connection $\nabla ,$ one has $d_{\nabla
}=d_{\gamma }-T,$ where 
\begin{equation*}
T:\Omega ^{1}\left( M\right) \rightarrow S^{2}\left( \Omega ^{1}(M\right)
)\subset \Omega ^{1}\left( M\right) \otimes \Omega ^{1}\left( M\right)
\end{equation*}%
is the \emph{affine deformation tensor }of the connection, and $S^{2}\left(
\Omega ^{1}(M\right) )$ is the module of the symmetric $(0,2)$-tensors on $M$%
.

In what follows, we shall use the notation $\nabla _{X}\left( \theta \right) 
\overset{\text{def}}{=}\left( d_{\nabla }\theta \right) \left( X\right) $
for the covariant derivative of a differential $1$-form $\theta $ along a
vector field $X$ with respect to the connection $\nabla .$

\begin{proposition}
\label{geodesics proposition}Let $d_{\nabla }:\Omega ^{1}\left( M\right)
\rightarrow \Omega ^{1}\left( M\right) \otimes \Omega ^{1}\left( M\right) $
be the covariant differential of a connection $\nabla $ in the cotangent
bundle of $M.$ Then the foliation $\left\{ \theta =0\right\} $ on $M$ given
by the differential $1$-form $\theta \in \Omega ^{1}\left( M\right) $
consists of geodesics of $\;\nabla $ if and only if 
\begin{equation*}
d_{\nabla }\left( \theta \right) =\alpha \otimes \theta +\theta \otimes \beta
\end{equation*}%
for some differential $1$-forms $\alpha ,\beta \in \Omega ^{1}\left(
M\right) .$
\end{proposition}

\begin{proof}
Let $\theta ^{\prime }$ be a differential $1$-form such that $\theta $ and $%
\theta ^{\prime }$ are linearly independent. Then 
\begin{equation*}
d_{\nabla }\left( \theta \right) =\alpha \otimes \theta +\theta \otimes
\beta +h\theta ^{\prime }\otimes \theta ^{\prime }.
\end{equation*}%
Assume that $X$ is a geodesic vector field on $M$ such that $\theta \left(
X\right) =0.$ Then $\nabla _{X}\left( \theta \right) $ must be equal to zero
on $X.$ But 
\begin{equation*}
d_{\nabla }\theta \left( X\right) =\beta \left( X\right) \theta +h\theta
^{\prime }\left( X\right) \theta ^{\prime }.
\end{equation*}%
Therefore, $h=0.$
\end{proof}

\begin{corollary}
The foliations $\left\{ \omega _{1}=0\right\} ,\left\{ \omega _{2}=0\right\}
,$ and $\left\{ \omega _{3}=0\right\} $ are geodesic with respect to the
Chern connection.
\end{corollary}

The problem of linearization of webs can be reformulated as follows: \emph{%
find a torsion-free flat connection such that the foliations of the web are
geodesic with respect to this connection.}

\begin{proposition}
Let $d_{\nabla }=d_{\gamma }-T:\Omega ^{1}\left( M\right) \rightarrow \Omega
^{1}\left( M\right) \otimes \Omega ^{1}\left( M\right) $ be the covariant
differential of a torsion-free connection $\nabla $ such that the foliations 
$\left\{ \omega _{p}=0\right\} ,\ p=1,2,3,$ are geodesic with respect to the
connection $\nabla $. Then%
\begin{eqnarray}
T &=&(T_{11}^{1}\omega _{1}\otimes \omega _{1}+T_{12}^{1}\left( \omega
_{1}\otimes \omega _{2}+\omega _{2}\otimes \omega _{1}\right) )\otimes
\partial _{1}  \notag \\
&&+(T_{22}^{2}\omega _{2}\otimes \omega _{2}+T_{12}^{2}\left( \omega
_{1}\otimes \omega _{2}+\omega _{2}\otimes \omega _{1}\right) )\otimes
\partial _{2},  \label{defTensorGeneral}
\end{eqnarray}%
where the components of the affine deformation tensor have the form 
\begin{equation}
T_{12}^{2}=\lambda _{1},\ \ T_{12}^{1}=\lambda _{2},\ T_{11}^{1}=2\lambda
_{1}+\mu ,\ T_{22}^{2}=2\lambda _{2}-\mu \   \label{Deformation tensor}
\end{equation}%
for some smooth functions $\lambda _{1},\lambda _{2},$ and $\mu .$
\end{proposition}

\begin{proof}
Due to Proposition \ref{geodesics proposition} and the requirement that the
foliations $\left\{ \omega _{1}=0\right\} \ $and $\left\{ \omega
_{2}=0\right\} \ $are geodesic, one gets (\ref{defTensorGeneral}). \ The
same requirement for the foliation $\left\{ \omega _{3}=0\right\} $ gives
the following relation for the components of the affine deformation tensor $%
T $: 
\begin{equation*}
T_{11}^{1}+T_{22}^{2}=2(T_{12}^{1}+T_{12}^{2}),
\end{equation*}%
and this implies (\ref{Deformation tensor}).
\end{proof}

Therefore, in order to linearize a $3$-web, one should find functions $%
\lambda _{1},\lambda _{2}\ \ $and $\mu $\ in such a way that the connection
corresponding to $d_{T}=d_{\gamma }-T,$\ where the affine deformation tensor 
$T$\ has form (\ref{Deformation tensor}), is flat.

The covariant differential $d_{T}$ has the following form: 
\begin{eqnarray*}
d_{T}\omega _{1} &=&-\omega _{1}\otimes \sigma _{11}-\omega _{2}\otimes
\sigma _{12}, \\
d_{T}\omega _{2} &=&-\omega _{1}\otimes \sigma _{21}-\omega _{2}\otimes
\sigma _{22},
\end{eqnarray*}%
where%
\begin{eqnarray*}
\sigma _{11} &=&\gamma +\left( 2\lambda _{1}+\mu \right) \omega _{1}+\lambda
_{2}\omega _{2}, \\
\ \sigma _{12} &=&\lambda _{2}\omega _{1}, \\
\sigma _{21} &=&\lambda _{1}\omega _{2}, \\
\ \sigma _{22} &=&\gamma +\lambda _{1}\omega _{1}+\left( 2\lambda _{2}-\mu
\right) \omega _{2}.
\end{eqnarray*}%
Using structure equations (\ref{web structure equations}), we get%
\begin{eqnarray*}
d_{T}^{2}\omega _{1} &=&\omega _{1}\otimes \left( \sigma _{21}\wedge \sigma
_{12}-d\sigma _{11}\right) +\omega _{2}\otimes \left( \sigma _{12}\wedge
\sigma _{11}+\sigma _{21}\wedge \sigma _{12}-d\sigma _{12}\right) , \\
d_{T}^{2}\omega _{2} &=&\omega _{1}\otimes \left( \sigma _{11}\wedge \sigma
_{21}+\sigma _{21}\wedge \sigma _{22}-d\sigma _{21}\right) +\omega
_{2}\otimes \left( \sigma _{12}\wedge \sigma _{21}-d\sigma _{22}\right) .
\end{eqnarray*}

Therefore, in order to obtain a flat torsion-free connection, components of
the affine deformation tensor must satisfy the following \emph{%
Akivis--Goldberg\ equations}: 
\begin{equation}  \label{AG1-equations}
\renewcommand{\arraystretch}{1.3} 
\begin{array}{ll}
d\sigma _{11}=\sigma _{21}\wedge \sigma _{12}, &  \\ 
d\sigma _{12}=\sigma _{12}\wedge \sigma _{11}+\sigma _{21}\wedge \sigma
_{12}, &  \\ 
d\sigma _{21}=\sigma _{11}\wedge \sigma _{21}+\sigma _{21}\wedge \sigma
_{22}, &  \\ 
d\sigma _{22}=\sigma _{12}\wedge \sigma _{21}. & 
\end{array}
\renewcommand{\arraystretch}{1}
\end{equation}

Because $\omega _{1}$\ and $\omega _{2}$\ are linearly independent,
equations (\ref{AG1-equations}) imply that 
\begin{equation}  \label{AG2-equations}
\renewcommand{\arraystretch}{1.3} 
\begin{array}{ll}
2\partial _{2}\left( \lambda _{1}\right) -\partial _{1}\left( \lambda
_{2}\right) +\partial _{2}\left( \mu \right) =K+\lambda _{1}\lambda
_{2}+\!g_{2}\left( 2\lambda _{1}+\mu \right) -g_{1}\!\lambda _{2}, &  \\ 
\partial _{2}\left( \lambda _{2}\right) =\lambda _{2}\left( g_{2}+\lambda
_{2}-\mu \right) ,\!\!\! &  \\ 
\partial _{1}\left( \lambda _{1}\right) =\lambda _{1}\left( g_{1}+\lambda
_{1}+\mu \right) ,\!\!\!\! &  \\ 
\partial _{2}\left( \lambda _{1}\right) -2\partial _{1}\left( \lambda
_{2}\right) +\partial _{1}\left( \mu \right) =K-\lambda _{1}\lambda
_{2}+\lambda _{1}g_{2}-g_{1}\left( 2\lambda _{2}-\mu \right) \!\!\!\!. & 
\end{array}
\renewcommand{\arraystretch}{1}
\end{equation}

\section{ Calculus of Covariant Derivatives}

Let $d_{\gamma }:\Omega ^{1}(M)\rightarrow \Omega ^{1}\left( M\right)
\otimes \Omega ^{1}\left( M\right) $ be the covariant differential with
respect to the Chern connection. It induces the connection $d_{\gamma
}^{\ast }:\mathcal{D}\left( M\right) \rightarrow \mathcal{D}\left( M\right)
\otimes \Omega ^{1}\left( M\right) $ in the tangent bundle, where 
\begin{eqnarray*}
d_{\gamma }^{\ast } &:&\partial _{1}\rightarrow \partial _{1}\otimes \gamma ,
\\
d_{\gamma }^{\ast } &:&\partial _{2}\rightarrow \partial _{2}\otimes \gamma .
\end{eqnarray*}

Denote by $\Theta ^{p,q}\left( M\right) =\left( \mathcal{D}\left( M\right)
\right) ^{\otimes p}\otimes \left( \Omega ^{1}\left( M\right) \right)
^{\otimes q}$ the module of tensors of type $\left( p,q\right) .\;$Then the
Chern connection induces the covariant differential 
\begin{equation*}
d_{\gamma }^{(p,q)}:\Theta ^{p,q}\left( M\right) \rightarrow \Theta
^{p+1,q}\left( M\right) ,
\end{equation*}%
where 
\begin{equation*}
d_{\gamma }^{(p,q)}:u\partial _{j_{1}}\otimes \cdots \otimes \partial
_{j_{p}}\otimes \omega _{i_{1}}\otimes \cdots \otimes \omega
_{i_{q}}\longmapsto \partial _{j_{1}}\otimes \cdots \otimes \partial
_{j_{p}}\otimes \omega _{i_{1}}\otimes \cdots \otimes \omega _{i_{q}}\otimes
\left( du+\left( p-q\right) \gamma u\right)
\end{equation*}%
and $u\in C^{\infty }\left( M\right) .$

We say that $u$ is of weight $q-p$ and call the form 
\begin{equation}
\delta ^{\left( p,q\right) }\left( u\right) \overset{\text{def}}{=}\delta
^{\left( q-p\right) }\left( u\right) =du-\left( q-p\right) u\gamma
\label{covariant differential Vadim}
\end{equation}%
the \emph{covariant differential }of $u.$

Decomposing the form $\delta ^{\left( q-p\right) }\left( u\right) $ in the
basis $\{\omega _{1},\omega _{2}\},$ we obtain 
\begin{equation*}
\delta ^{\left( q-p\right) }\left( u\right) =\delta _{1}^{\left( q-p\right)
}\left( u\right) ~\omega _{1}+\delta _{2}^{\left( q-p\right) }\left(
u\right) ~\omega _{2},
\end{equation*}%
where 
\begin{equation}  \label{covariant derivatives weight k}
\renewcommand{\arraystretch}{1.5} 
\begin{array}{ll}
\delta _{1}^{\left( q-p\right) }\left( u\right)
\!\!\!\!~~=~~\!\!\!\!\partial _{1}\left( u\right) -\left( q-p\right) g_{1}u,
&  \\ 
\delta _{2}^{\left( q-p\right) }\left( u\right)
\!\!\!\!~~=~\!\!\!\!~\partial _{2}\left( u\right) -\left( q-p\right) g_{2}u
& 
\end{array}
\renewcommand{\arraystretch}{1}
\end{equation}%
are the covariant derivatives of $u$ with respect to the Chern connection.

Note that $\delta _{1}^{\left( q-p\right) }\left( u\right) $ and $\delta
_{2}^{\left( q-p\right) }\left( u\right) $ are of weight $q-p+1.$

\begin{lemma}
For any $s=0,\pm 1,\pm 2,...,$ the relation 
\begin{equation}
\delta _{2}^{\left( s+1\right) }\circ \delta _{1}^{\left( s\right) }-\delta
_{1}^{\left( s+1\right) }\circ \delta _{2}^{\left( s\right) }=sK
\label{covariant commutator}
\end{equation}%
holds for the commutator.
\end{lemma}

\begin{proof}
We have 
\begin{equation*}
\delta _{2}^{\left( s+1\right) }\circ \delta _{1}^{\left( s\right)
}=\partial _{2}\partial _{1}-sg_{1}\partial _{2}-\left( s+1\right)
g_{2}\partial _{1}+s\left( s+1\right) g_{1}g_{2}-s\partial _{2}(g_{1})
\end{equation*}%
and 
\begin{equation*}
\delta _{1}^{\left( s+1\right) }\circ \delta _{2}^{\left( s\right)
}=\partial _{1}\partial _{2}-sg_{2}\partial _{1}-\left( s+1\right)
g_{1}\partial _{2}+s\left( s+1\right) g_{1}g_{2}-s\partial _{1}(g_{2}).
\end{equation*}%
The statement follows now from (\ref{curvature main}).
\end{proof}

Note that the curvature $K$ is of weight two, while $\lambda _{1},\lambda
_{2}$ and $\mu $ are of weight one.

The classical Leibnitz rule leads to the corresponding rule for weighted
functions.

\begin{lemma}[Leibnitz rule]
Let $u$ be of weight $k$ and $v$ be of weight $l$. Then%
\begin{equation*}
\delta _{i}^{\left( k+l\right) }\left( uv\right) =\delta _{i}^{\left(
k\right) }\left( u\right) ~v+u~\delta _{i}^{\left( l\right) }\left( v\right)
.
\end{equation*}
\end{lemma}

In what follows, we shall omit the superscript indicating the weight in the
cases when the weight is known. For example, we shall write $\delta _{1}K$
instead of $\delta _{1}^{\left( 2\right) }K,$ or $\delta _{1}\mu $ instead
of $\delta _{1}^{\left( 1\right) }\mu $.

\section{Differential Invariants and Rigidity of 3-Webs}

As we have noted above, the curvature $K$ is a relative invariant of weight
two of a $3$-web $W.$ The covariant derivatives of $K$ are relative
invariants of weight three. The invariants (\ref{a-invariants}) can be
written in terms of the curvature $K$\ as follows: 
\begin{equation*}
a_{1}=\frac{-\delta _{1}K}{2K^{\frac{3}{2}}},\ \ \ \ \ \ a_{2}=-\frac{\delta
_{2}K}{2K^{\frac{3}{2}}}.
\end{equation*}%
They are absolute invariants of a $3$-web $W$ with nonvanishing curvature $K$%
.

Hence all the derivatives 
\begin{equation*}
a_{1}^{i,j}=\nabla _{1}^{i}\nabla _{2}^{j}(a_{1})\ \ \ \text{and\ \ \ }%
a_{2}^{i,j}=\nabla _{1}^{i}\nabla _{2}^{j}(a_{2})\text{\ }
\end{equation*}%
are absolute invariants too; here $i,j=0,1,2,...$

It is easy to see that they are differential operators with respect to the
web function $f$ of order $i+j+4.$

Note also that condition (\ref{normalized K}), 
\begin{equation*}
\nabla _{1}\left( a_{2}\right) -\nabla _{2}\left( a_{1}\right) =1,
\end{equation*}%
gives the differential relations between the invariants $a_{1}^{i,j}\ $and$\
a_{2}^{i,j}$.

In particular, it follows that there are no $3$-webs with constant
invariants $a_{1}$ and $a_{2}.$

The following theorem is valid (cf. \cite{B 55}, \S 13 and \cite{BB 38}, \S
20).

\begin{theorem}
The differential invariants $a_{1}^{i,j}\ $and$\ a_{2}^{i,j}$ form a
complete system of differential invariants of $\;3$-webs with nonvanishing
curvature, that is, any differential invariant of such $3$-webs is a
function of a finite number of invariants from the system $\left\{
a_{1}^{i,j},a_{2}^{i,j}\right\} ,\ i,j=0,1,2,...$
\end{theorem}

We say that a $3$-web $W$ is \emph{locally rigid in a domain }$D\subset M$
if for any two distinct points \ $p,q\in D$ there is no local diffeomorphism 
$\phi $ sending $p$ to $q$ and transforming the web $W$ in a neighborhood of 
$p$ into the web $W$ in a neighborhood of $q.$

The problem of local rigidity can be viewed as a generalized Gronwall
conjecture (see the description of the Gronwall conjecture for linearizable
webs in Section 8 or in \cite{B 55}, \S 17).

It is easy to see that locally rigid webs do not have nontrivial
(infinitesimal) automorphisms.

Let $W$ be\ a $3$-web defined in some neighborhood $D$ of the point $p,$ let 
$\theta _{1},\theta _{2}$ and $\alpha $ be its invariant differential $1$%
-forms, and let $a_{1},a_{2}$ be its absolute differential invariants.
Denote by $\overline{W}$ a copy of $W$ with corresponding forms $\overline{%
\theta }_{1},\overline{\theta }_{2}$ , $\overline{\alpha }$ and invariants $%
\overline{a}_{1},\overline{a}_{2}.$

On the product $D\times D,$ we consider the $1$-forms%
\begin{equation*}
\Theta _{1}=\overline{\theta }_{1}-\theta _{1},\ \Theta _{2}=\overline{%
\theta }_{2}-\theta _{2},\ \ \aleph =\overline{\alpha }-\alpha
\end{equation*}%
and the functions 
\begin{equation*}
A_{1}=\overline{a_{1}}-a_{1},\ A_{2}=\overline{a_{2}}-a_{2}.
\end{equation*}%
Then the graph $G_{\phi }\subset D\times D$ of a local diffeomorphism $\phi
:D\rightarrow D,$ $\phi \left( p\right) =q,$ transforming $W$ in a
neighborhood of $p$ into $W$ in a neighborhood of $q$ is an integral surface
of the differential system 
\begin{equation}
\Theta _{1}=0,\ \Theta _{2}=0,\ \ \aleph =0  \label{equivalence dif. system}
\end{equation}%
such that 
\begin{equation}
\left. A_{1}\right| _{G_{\phi }}=0,\ \left. A_{2}\right| _{G_{\phi }}=0.
\label{function equivalence condition}
\end{equation}%
Assume that the functions $a_{1}$ and $a_{2}$ are functionally independent
in $D,$ and $D$ is sufficiently small. Then the invariants $a_{1}$ and $%
a_{2} $ can be viewed as coordinates on $D,$ and therefore the distinct
points $p$ and $q$ have distinct coordinates.$\ $ This means that the web $W$
is locally rigid.

Let us assume that there is a functional dependence between the invariants $%
a_{1}$ and $a_{2}$, say, $a_{2}=F\left( a_{1}\right) .$ Then (\ref{function
equivalence condition}) determines a $3$-dimensional manifold $N$ such that
the graphs $G_{\phi }$ are integral surfaces of differential system (\ref%
{equivalence dif. system}) on $N.$

For the system 
\begin{equation*}
\left. \Theta _{1}\right| _{N}=0,\left. \Theta _{2}\right| _{N}=0,\ \ \left.
\aleph \right| _{N}=0
\end{equation*}%
to have two-dimensional integral manifolds, it is necessary and sufficient
that the forms $\left. \Theta _{1}\right| _{N},\left. \Theta _{2}\right|
_{N} $ and $\left. \aleph \right| _{N}$ are proportional. In fact, the
distribution defined by the above system should be two-dimensional and
completely integrable. This follows from the fact that proportionality of
these forms implies complete integrability of the system.

Indeed, let $\left. \Theta _{1}\right| _{N}\wedge \left. \Theta _{2}\right|
_{N}=0.$ Then 
\begin{equation*}
\left. \aleph \right| _{N}=a_{1}\left. \Theta _{1}\right| _{N}+a_{2}\left.
\Theta _{2}\right| _{N},
\end{equation*}%
and therefore $\left. \Theta _{1}\right| _{N}\wedge \left. \aleph \right|
_{N}=\left. \Theta _{2}\right| _{N}\wedge \left. \aleph \right| _{N}=0.$

Moreover,%
\begin{equation*}
d\left. \Theta _{i}\right| _{N}=\left. \Theta _{i}\right| _{N}\wedge \left. 
\overline{\alpha }\right| _{N}+\left. \theta _{i}\right| _{N}\wedge \left.
\aleph \right| _{N},
\end{equation*}%
and hence the system is completely integrable.

Summarizing, we arrive at the following theorem.

\begin{theorem}

\begin{enumerate}
\item[$(i)$] Let $W$ be a $3$-web defined in a domain $D$ in which the
invariants $a_{1}$ and $a_{2}$ are functionally independent and form a
coordinate system. Then $W$ is locally rigid in $D$.
\end{enumerate}

\begin{itemize}
\item[$(ii)$] Let the invariants $a_{1}\;$and $a_{2}$ be functionally
dependent in some domain $D,$ say, $a_{2}=F\left( a_{1}\right) ,$ for a
smooth function $F,$ but the differential $3$-form 
\begin{equation}
\Theta _{1}\wedge \Theta _{2}\wedge dA_{1}\neq 0
\end{equation}%
at points of the manifold $\{\left. \left( p,q\right) \right| \ a_{1}\left(
p\right) =\overline{a_{1}}\left( q\right) ,\ p\neq q\}\subset D\times D.$
Then $W$ is locally rigid in this domain.
\end{itemize}
\end{theorem}

We say that a vector field $X$ is an \emph{infinitesimal automorphism} of a $%
3$-web $W$ if the one-parameter group of shifts along $X$ consists of
diffeomorphisms preserving $W.$ A $3$-web $W$ is said to be \emph{%
infinitesimally rigid }if $W$ has the trivial infinitesimal automorphism ($%
X=0$) only.

In terms of the invariant forms $\theta _{1}$ and $\theta _{2},$ this means
that the following Lie equations 
\begin{equation*}
L_{X}\left( \theta _{1}\right) =0,\ L_{X}\left( \theta _{1}\right) =0
\end{equation*}%
hold. Here $L_{X}$ is the Lie derivative along $X.$

Let 
\begin{equation*}
X=X_{1}\nabla _{1}+X_{2}\nabla _{2}
\end{equation*}%
be the decomposition of $X$ in the basis $\left\{ \nabla _{1},\nabla
_{2}\right\} .$ Using structure equations (\ref{standard structure equations}%
), one can rewrite the Lie equations as follows:%
\begin{eqnarray*}
dX_{1} &=&a_{2}X_{2}\theta _{1}-a_{2}X_{1}\theta _{2}, \\
dX_{2} &=&-a_{1}X_{2}\theta _{1}+a_{1}X_{1}\theta _{2},
\end{eqnarray*}%
or%
\begin{equation}  \label{Lie component equations}
\renewcommand{\arraystretch}{1.3} 
\begin{array}{ll}
\nabla _{1}\left( X_{1}\right) =a_{2}X_{2},\ \ \ \ \nabla _{2}\left(
X_{1}\right) =-a_{2}X_{1}, &  \\ 
\nabla _{1}\left( X_{2}\right) =-a_{1}X_{2},\ \ \nabla _{2}\left(
X_{2}\right) =a_{1}X_{1}. & 
\end{array}
\renewcommand{\arraystretch}{1}
\end{equation}%
The compatibility conditions for these equations follow from (\ref{nabla
commutator}). Namely, applying the operators from the left- and right-hand
sides of (\ref{nabla commutator}) to $X_{1}$ and $X_{2},$ we get%
\begin{equation}  \label{Symmetry system}
\renewcommand{\arraystretch}{1.3} 
\begin{array}{ll}
\nabla _{1}\left( a_{2}\right) X_{1}+\nabla _{2}\left( a_{2}\right) X_{2}=0,
&  \\ 
\nabla _{1}\left( a_{1}\right) X_{1}+\nabla _{2}\left( a_{1}\right) X_{2}=0.
& 
\end{array}
\renewcommand{\arraystretch}{1}
\end{equation}%
This implies the following theorem.

\begin{theorem}[Infinitesimal Rigidity of 3-Webs]
Let $W$ be a $3$-web given in a domain $D,$ and let the invariant%
\begin{equation*}
J=\det 
\begin{Vmatrix}
\nabla _{1}\left( a_{1}\right) & \nabla _{1}\left( a_{2}\right) \\ 
\nabla _{2}\left( a_{1}\right) & \nabla _{2}\left( a_{2}\right)%
\end{Vmatrix}%
\end{equation*}%
be nonvanishing in $D.$ Then $W$ is infinitesimally rigid in $D.$
\end{theorem}

Let us assume now that $J$ identically equals zero in $D.$ As we have seen
earlier, the entries of the above matrix do not vanish simultaneously, that
is, the rank of the matrix equals one.

Hence system (\ref{Symmetry system}) has solutions of the form 
\begin{equation*}
X=s\left( \nabla _{2}\left( a_{2}\right) ~\nabla _{1}-\nabla _{1}\left(
a_{2}\right) ~\nabla _{2}\right)
\end{equation*}%
for some smooth function $s.$

Substituting this expression into system (\ref{Lie component equations}), we
get%
\begin{equation}  \label{sym s system}
\renewcommand{\arraystretch}{2.2} 
\begin{array}{ll}
\nabla _{1}\left( s\right) =-\displaystyle\frac{a_{2}\nabla _{1}\left(
a_{2}\right) +\nabla _{1}\nabla _{2}\left( a_{2}\right) }{\nabla _{2}\left(
a_{2}\right) }s, &  \\ 
\nabla _{1}\left( s\right) =-\displaystyle\frac{a_{1}\nabla _{1}\left(
a_{2}\right) +\nabla _{1}^{2}\left( a_{2}\right) }{\nabla _{1}\left(
a_{2}\right) }s, &  \\ 
\nabla _{2}\left( s\right) =-\displaystyle\frac{a_{2}\nabla _{2}\left(
a_{2}\right) +\nabla _{2}^{2}\left( a_{2}\right) }{\nabla _{2}\left(
a_{2}\right) }s, &  \\ 
\nabla _{2}\left( s\right) =-\displaystyle\frac{a_{1}\nabla _{2}\left(
a_{2}\right) +\nabla _{2}\nabla _{1}\left( a_{2}\right) }{\nabla _{1}\left(
a_{2}\right) }s. & 
\end{array}
\renewcommand{\arraystretch}{1}
\end{equation}%
It follows that 
\begin{equation}  \label{1st inv relations}
\renewcommand{\arraystretch}{1.3} 
\begin{array}{ll}
a_{2}(\nabla _{1}\left( a_{2}\right) )^{2}+\nabla _{1}\nabla _{2}\left(
a_{2}\right) ~\nabla _{1}\left( a_{2}\right) =a_{1}\nabla _{1}\left(
a_{2}\right) ~\nabla _{2}\left( a_{2}\right) +\nabla _{1}^{2}\left(
a_{2}\right) ~\nabla _{2}\left( a_{2}\right) , &  \\ 
a_{2}\nabla _{2}\left( a_{2}\right) \nabla _{1}\left( a_{2}\right) +\nabla
_{2}^{2}\left( a_{2}\right) \nabla _{1}\left( a_{2}\right) =a_{1}(\nabla
_{2}\left( a_{2}\right) )^{2}+\nabla _{2}\nabla _{1}\left( a_{2}\right)
\nabla _{2}\left( a_{2}\right) . & 
\end{array}
\renewcommand{\arraystretch}{1}
\end{equation}

The compatibility conditions for the above system take the form:%
\begin{eqnarray*}
&&\nabla _{2}\left( \frac{a_{2}\nabla _{1}\left( a_{2}\right) +\nabla
_{1}\nabla _{2}\left( a_{2}\right) }{\nabla _{2}\left( a_{2}\right) }\right)
-\nabla _{1}\left( \frac{a_{2}\nabla _{2}\left( a_{2}\right) +\nabla
_{2}^{2}\left( a_{2}\right) }{\nabla _{2}\left( a_{2}\right) }\right) \\
&=&-a_{2}\frac{a_{2}\nabla _{1}\left( a_{2}\right) +\nabla _{1}\nabla
_{2}\left( a_{2}\right) }{\nabla _{2}\left( a_{2}\right) }+a_{1}\frac{%
a_{2}\nabla _{2}\left( a_{2}\right) +\nabla _{2}^{2}\left( a_{2}\right) }{%
\nabla _{2}\left( a_{2}\right) }
\end{eqnarray*}%
or 
\begin{equation*}
~\nabla _{2}\nabla _{1}\nabla _{2}\left( a_{2}\right) +a_{2}\nabla
_{1}\nabla _{2}\left( a_{2}\right) =a_{1}\nabla _{2}^{2}\left( a_{2}\right)
+a_{2}\nabla _{1}\left( a_{2}\right) \nabla _{2}\left( a_{2}\right) +\nabla
_{1}\nabla _{2}^{2}\left( a_{2}\right) .
\end{equation*}

\begin{theorem}
Let $W$ be a $3$-web such that $J=0$, and suppose that the invariants $a_{1}$
and $a_{2}$ satisfy the relations%
\begin{eqnarray*}
a_{2}(\nabla _{1}\left( a_{2}\right) )^{2}+\nabla _{1}\nabla _{2}\left(
a_{2}\right) ~\nabla _{1}\left( a_{2}\right) &=&a_{1}\nabla _{1}\left(
a_{2}\right) ~\nabla _{2}\left( a_{2}\right) +\nabla _{1}^{2}\left(
a_{2}\right) ~\nabla _{2}\left( a_{2}\right) , \\
a_{2}\nabla _{2}\left( a_{2}\right) \nabla _{1}\left( a_{2}\right) +\nabla
_{2}^{2}\left( a_{2}\right) \nabla _{1}\left( a_{2}\right) &=&a_{1}(\nabla
_{2}\left( a_{2}\right) )^{2}+\nabla _{2}\nabla _{1}\left( a_{2}\right)
\nabla _{2}\left( a_{2}\right) , \\
\nabla _{2}\nabla _{1}\nabla _{2}\left( a_{2}\right) +a_{2}\nabla _{1}\nabla
_{2}\left( a_{2}\right) &=&a_{1}\nabla _{2}^{2}\left( a_{2}\right)
+a_{2}\nabla _{1}\left( a_{2}\right) \nabla _{2}\left( a_{2}\right) +\nabla
_{1}\nabla _{2}^{2}\left( a_{2}\right) .
\end{eqnarray*}%
Then there is a nontrivial infinitesimal automorphism of $W$ which is unique
up to a factor and has the form 
\begin{equation*}
X=s\left( \nabla _{2}\left( a_{2}\right) ~\nabla _{1}-\nabla _{1}\left(
a_{2}\right) ~\nabla _{2}\right) ,
\end{equation*}%
where the function $s$ is a solution of $(\ref{sym s system})$.
\end{theorem}

\subsection{Examples}

\textbf{Example 4 }Consider the $3$-web $W$\ given by the web function

\begin{equation*}
f=x+\sqrt{x^{2}-y}
\end{equation*}%
in the domain $\left\{ x>0,\;y>0,\;y<x^{2}\right\} \;$(cf. Example $3$).$\;$

As we saw in Example $3$, this web is generated by two families of
coordinate lines $\left\{ x=\limfunc{const}\right\} $, $\left\{ y=\limfunc{%
const}\right\} $ and the tangents to the parabola $y=x^{2}.$

For this web, we have 
\begin{eqnarray*}
\omega _{1} &=&-\frac{f~dx}{f-x},\ \ \omega _{2}=\frac{dy}{2\left(
f-x\right) },\ \gamma =\frac{x\left( -2fdx+dy\right) }{2\left( f-x\right)
\left( y-xf\right) }, \\
H &=&\frac{x}{y-xf},\ \ K=\frac{2x^{2}f-y\left( f+x\right) }{f\left(
xf-y\right) ^{2}}, \\
\theta _{1} &=&-\frac{\sqrt{f}dx}{f-x},\ \theta _{2}=\frac{dy}{2\sqrt{f}%
\left( f-x\right) },\ \alpha =\frac{\left( f+2x\right) ~dx}{2\left(
f-x\right) ^{3/2}}-\frac{\left( 2f+x\right) dy}{4\sqrt{f\left( f-x\right) }},
\\
a_{1} &=&-\frac{f+2x}{2\sqrt{f\left( f-x\right) }},\ \ a_{2}=-\frac{2f+x}{2%
\sqrt{f\left( f-x\right) }}
\end{eqnarray*}%
Note that $da_{1}\wedge da_{2}=0.$ Hence the invariants $a_{1}$ and $a_{2}$
are functionally dependent. The dependence is%
\begin{equation*}
8a_{1}^{2}-5a_{2}^{2}+4a_{1}\left( a_{2}^{2}-1\right) \sqrt{a_{1}^{2}+6}%
+a_{2}\left( 4a_{1}^{2}-1\right) \sqrt{a_{2}^{2}+3}+3=0.
\end{equation*}%
Conditions (\ref{function equivalence condition}) mean that 
\begin{equation*}
a_{1}\left( x,y\right) =a_{1}\left( \overline{x},\overline{y}\right)
\end{equation*}%
or 
\begin{equation*}
\frac{y}{x^{2}}=\frac{\overline{y}}{\overline{x}^{2}}.
\end{equation*}%
Then 
\begin{equation*}
\Theta _{1}=\frac{\sqrt{f}~(\sqrt{x}d\overline{x}-\sqrt{\overline{x}}dx)}{%
(f-x)\sqrt{\overline{x}}},
\end{equation*}%
and%
\begin{equation*}
\Theta _{2}=\frac{x\sqrt{x}}{2\sqrt{f}\left( f-x\right) }\Biggl[\frac{\sqrt{%
\overline{x}}-\sqrt{x}~}{x^{2}}dy+\frac{2y}{x^{3}\sqrt{\overline{x}}}(xd%
\overline{x}-\overline{x}dx)\Biggr].
\end{equation*}

It is easy to check that on the manifold $N,$ the condition $\Theta
_{1}\wedge \Theta _{2}=0$ holds if and only if $x=\overline{x}$ and
consequently $y=\overline{y}.$

In other words, this web is locally rigid.

\bigskip

\textbf{Example 5} Consider the $3$-web $W$\ given by the web function

\begin{equation}
f=(x+y)e^{-x}.  \label{Exp}
\end{equation}%
$\;$

This web is generated by two families of coordinate lines $\left\{ x=%
\limfunc{const}\right\} $, $\left\{ y=\limfunc{const}\right\} $ and the
level sets of the function $f.$

Let $t=1-x-y.$ Then for web (\ref{Exp}) one has 
\begin{eqnarray*}
\omega _{1} &=&-te^{-x}dx,\ \ \omega _{2}=-e^{-x}dy,\ \gamma =dx+\frac{dy}{t}%
, \\
H &=&-\frac{e^{x}}{t},\ \ K=\frac{e^{2x}}{t^{3}}, \\
\theta _{1} &=&-\frac{dx}{\sqrt{t}},\ \theta _{2}=-\frac{dy}{t^{3/2}},\
\alpha =-\frac{3dx+dy}{2t}, \\
a_{1} &=&\frac{3}{2\sqrt{t}},\ \ a_{2}=\frac{\sqrt{t}}{2}.
\end{eqnarray*}%
Note that $da_{1}\wedge da_{2}=0.$ Hence the invariants $a_{1}$ and $a_{2}$
are functionally dependent:%
\begin{equation*}
a_{1}a_{2}=\frac{3}{4}.
\end{equation*}%
The three-dimensional manifold $N$ is defined by 
\begin{equation*}
x+y=\overline{x}+\overline{y},
\end{equation*}%
and the differential $1$-forms are%
\begin{equation*}
\left. \Theta _{1}\right| _{N}=\frac{dx-d\overline{x}}{\sqrt{t}},\left.
\Theta _{2}\right| _{N}=-\frac{dx-d\overline{x}}{t^{3/2}},\ \ \left. \aleph
\right| _{N}=\frac{dx-d\overline{x}}{t}.
\end{equation*}%
Therefore the integral surfaces are given by the equations:%
\begin{equation*}
\overline{x}=x+c,\ \overline{y}=y-c,
\end{equation*}%
and the requirement $\phi \left( p\right) =p$ implies $c=0.$

Therefore web (\ref{Exp}) is not locally rigid. Note that the vector field 
\begin{equation*}
X=\frac{\partial }{\partial x}-\frac{\partial }{\partial y}
\end{equation*}%
is the infinitesimal symmetry of web (\ref{Exp}).

\section{Akivis--Goldberg Equations}

Using the covariant derivatives instead of the partial derivatives, we write
equations (\ref{AG2-equations}) as follows:%
\begin{eqnarray*}
2\delta _{2}\left( \lambda _{1}\right) -\delta _{1}\left( \lambda
_{2}\right) +\delta _{2}\left( \mu \right) &=&K+\lambda _{1}\lambda _{2}, \\
\delta _{2}\left( \lambda _{2}\right) &=&\lambda _{2}\left( \lambda _{2}-\mu
\right) , \\
\delta _{1}\left( \lambda _{1}\right) &=&\lambda _{1}\left( \lambda _{1}+\mu
\right) , \\
\delta _{2}\left( \lambda _{1}\right) -2\delta _{1}\left( \lambda
_{2}\right) +\delta _{1}\left( \mu \right) &=&K-\lambda _{1}\lambda _{2}.
\end{eqnarray*}

Solving this system with respect to the covariant derivatives of $\lambda
_{1}$ and $\lambda _{2}$, we obtain the following system of PDEs: 
\begin{eqnarray*}
\delta _{1}\left( \lambda _{1}\right) &&\!\!\!\!=~~\!\!\!\!\lambda
_{1}\left( \lambda _{1}+\mu \right) , \\
\delta _{2}\left( \lambda _{1}\right) &&\!\!\!\!=\lambda _{1}\lambda
_{2}+\!\!\!\!\;\;\frac{K}{3}+\frac{1}{3}\delta _{1}\left( \mu \right) -\frac{%
2}{3}\delta _{2}\left( \mu \right) , \\
\delta _{1}\left( \lambda _{2}\right) &&\!\!\!\!=\!\!\!\!~~\lambda
_{1}\lambda _{2}-\frac{K}{3}+\frac{2}{3}\delta _{1}\left( \mu \right) -\frac{%
1}{3}\delta _{2}\left( \mu \right) , \\
\delta _{2}\left( \lambda _{2}\right) &&\!\!\!\!=\!\!\!\!~~\lambda
_{2}\left( \lambda _{2}-\mu \right) .
\end{eqnarray*}

We shall look at the above system as a system of partial differential
equations with respect to the functions $\lambda _{1}$ and $\lambda _{2}$
provided that $\mu $ is given.

From(\ref{covariant commutator}) we get the compatibility conditions for
this system: 
\begin{equation*}
\delta _{1}(\delta _{2}\left( \lambda _{i}\right) )-\delta _{2}(\delta
_{1}\left( \lambda _{i}\right) )+K\lambda _{i}=0,
\end{equation*}%
where $i=1,2.$

After a series of straightforward computations, we obtain the following two
compatibility equations: 
\begin{equation}
I_{1}\left( \mu \right) =0,\ I_{2}\left( \mu \right) =0,\ 
\label{Compatibilty Equations}
\end{equation}%
where $I_{1}(\mu )\;$and $I_{1}(\mu )\;$have the form%
\begin{equation*}
I_{1}(\mu )=\delta _{1}^{2}(\mu )-2\delta _{1}\delta _{2}\left( \mu \right)
-\mu \delta _{1}\left( \mu \right) +2\mu \delta _{2}\left( \mu \right) -\mu
K+\delta _{1}(K)
\end{equation*}%
and%
\begin{equation*}
I_{2}(\mu )=\delta _{2}^{2}(\mu )-2\delta _{1}\delta _{2}\left( \mu \right)
-2\mu \delta _{1}\left( \mu \right) +\mu \delta _{2}\left( \mu \right) -\mu
K+\delta _{2}(K).
\end{equation*}

We shall use the symmetrized derivatives. Namely, let%
\begin{equation*}
\delta _{ij}=\frac{1}{2}(\delta _{i}\delta _{j}+\delta _{j}\delta _{i})
\end{equation*}%
be the symmetrized mixed second derivatives.

Then for functions of weight one, we have%
\begin{eqnarray*}
\delta _{12} &=&\delta _{1}\delta _{2}+\frac{K}{2}, \\
\delta _{21} &=&\delta _{1}\delta _{2}-\frac{K}{2},
\end{eqnarray*}%
and the expressions for $I_{1}(\mu )\;$and $I_{1}(\mu )\;$can be written as
follows:%
\begin{equation}  \label{mu-invariants}
\renewcommand{\arraystretch}{1.3} 
\begin{array}{lll}
I_{1}\left( \mu \right) & = & \delta _{11}\left( \mu \right) -2\delta
_{12}\left( \mu \right) -\mu \delta _{1}\left( \mu \right) +2\mu \delta
_{2}\left( \mu \right) +\delta _{1}\left( K\right) , \\ 
I_{2}\left( \mu \right) & = & \delta _{22}\left( \mu \right) -2\delta
_{12}\left( \mu \right) -2\mu \delta _{1}\left( \mu \right) +\mu \delta
_{2}\left( \mu \right) +\delta _{2}\left( K\right) .%
\end{array}
\renewcommand{\arraystretch}{1}
\end{equation}

We summarize these results in the following theorem.

\begin{theorem}[\protect\cite{AGL}]
The Akivis--Goldberg equations as differential equations with respect to the
components $T_{12}^{1}=\lambda _{2}$ and $T_{12}^{2}=\lambda _{1}$ of the
affine deformation tensor $T$\ are compatible if and only if the function $%
\mu $ satisfies the following differential equations: 
\begin{equation}
I_{1}\left( \mu \right) =0,\ I_{2}\left( \mu \right) =0.
\label{AG-compatibility eq}
\end{equation}%
If conditions $(\ref{AG-compatibility eq})$ are valid, then system $(\ref%
{AG1-equations})$ of \ PDEs is a Frobenius-type system, and for given values 
$\lambda _{1}\left( x_{0}\right) $ and $\lambda _{2}\left( x_{0}\right) $ \
at a point $x_{0}\in M^{2},$ there is $($a unique$)$ smooth solution of the
system in some neighborhood of $x_{0}.$
\end{theorem}

Let us denote by $\tau $ the following involution:%
\begin{equation*}
\tau :\left( x,y,\mu ,K\right) \rightarrow \left( y,x,-\mu ,-K\right) .
\end{equation*}%
Then one can check that%
\begin{equation*}
\tau \left( I_{1}\right) =I_{2}.
\end{equation*}

\section{ Calculus in Jet Spaces of Weighted Functions}

\subsection{Cartan's Forms in Nonholonomic Coordinates}

Let $\mathbb{J}^{r}\left( s\right) $ be the space of $r$-jets of weight $s$
functions in the plane $\mathbb{R}^{2}.$ We shall use the coordinates $%
\left( x,y,u,p_{1},p_{2},...,p_{i_{1}...i_{l}},...\right) $ in this space
corresponding to the symmetrized covariant derivatives, that is,%
\begin{eqnarray*}
u(j_{r}\left( h\right) ) &=&h,\ p_{1}\left( j_{r}\left( h\right) \right)
=\delta _{1}\left( h\right) ,\ p_{2}\left( j_{r}\left( h\right) \right)
=\delta _{2}\left( h\right) , \\
&&p_{i_{1}...i_{l}}\left( j_{r}\left( h\right) \right) =\delta
_{i_{1}...i_{l}}\left( h\right) ,\;...
\end{eqnarray*}%
Here $j_{r}\left( h\right) $ is the $r$-jet of the function $h.$ The
function $u\ $is of weight $s,$ and $\delta _{i_{1}...i_{l}}$ is its
symmetrized covariant derivative of order $i_{1}+\cdots +i_{l}.$

In what follows, we shall denote the symmetrized covariant derivatives of
the curvature function $K$\ by%
\begin{equation*}
K_{i_{1}....i_{l}}\overset{\text{def}}{=}\delta _{i_{1}...i_{l}}\left(
K\right) .
\end{equation*}

We describe now the Cartan distribution (see \cite{KLV} or \cite{ALV}) in $%
\mathbb{J}^{r}\left( s\right) $ in these coordinates. Let us begin with $%
\mathbb{J}^{1}\left( s\right) .$ The formula%
\begin{equation*}
df=\left( \delta _{1}f+sg_{1}f\right) \omega _{1}+\left( \delta
_{2}f+sg_{2}f\right) \omega _{2},
\end{equation*}%
where $f$ is a function of weight $s$, shows that the contact form on $%
\mathbb{J}^{1}\left( s\right) $ can be expressed as%
\begin{eqnarray*}
\varepsilon _{0} &=&du-\left( p_{1}+sg_{1}u\right) \omega _{1}-\left(
p_{2}+sg_{2}u\right) \omega _{2} \\
&=&du-su\gamma -p_{1}\omega _{1}-p_{2}\omega _{2}.
\end{eqnarray*}%
To find the Cartan forms on $\mathbb{J}^{2}\left( s\right) ,$ we shall use
the relations 
\begin{eqnarray*}
&&\delta _{1}\delta _{2}-\delta _{2}\delta _{1}=-wK, \\
&&\delta _{12}=\frac{1}{2}(\delta _{1}\delta _{2}+\delta _{2}\delta _{1}),
\end{eqnarray*}%
which hold for functions of weight $w.$

These formulae imply that%
\begin{equation}  \label{Composition deltas}
\renewcommand{\arraystretch}{1.3} 
\begin{array}{lll}
\delta _{1}\delta _{2} & = & \delta _{12}-\frac{1}{2}wK, \\ 
\delta _{2}\delta _{1} & = & \delta _{12}+\frac{1}{2}wK%
\end{array}
\renewcommand{\arraystretch}{1}
\end{equation}%
and give the following representation of the second-order Cartan forms:%
\begin{eqnarray*}
\varepsilon _{1} &=&dp_{1}-\left( s+1\right) p_{1}\gamma -p_{11}\omega
_{1}-(p_{12}+\frac{1}{2}sKu)\omega _{2}, \\
\varepsilon _{2} &=&dp_{2}-\left( s+1\right) p_{2}\gamma -(p_{12}-\frac{1}{2}%
sKu)~\omega _{1}-p_{22}\omega _{2}.
\end{eqnarray*}

To obtain the Cartan forms on the next jet space $\mathbb{J}^{3}\left(
s\right) ,$ we need the following relations:%
\begin{equation}  \label{delta 1-2 products}
\renewcommand{\arraystretch}{1.5} 
\begin{array}{lll}
\delta _{1}\delta _{12} & = & \delta _{112}-\frac{1}{6}\left( 3s+2\right)
K\delta _{1}-\frac{1}{6}sK_{1}, \\ 
\delta _{2}\delta _{12} & = & \delta _{122}+\frac{1}{6}\left( 3s+2\right)
K\delta _{2}+\frac{1}{6}sK_{2}, \\ 
\delta _{1}\delta _{22} & = & \delta _{122}-\frac{1}{3}\left( 3s+2\right)
K\delta _{2}-\frac{1}{3}sK_{2}, \\ 
\delta _{2}\delta _{11} & = & \delta _{112}+\frac{1}{3}\left( 3s+2\right)
K\delta _{1}+\frac{1}{3}sK_{1},%
\end{array}
\renewcommand{\arraystretch}{1}
\end{equation}%
which follow from (\ref{Composition deltas}).

These relations allow us to represent the third-order Cartan forms:%
\begin{eqnarray*}
\varepsilon _{11} &=&dp_{11}-\left( s+2\right) p_{11}~\gamma -p_{111}\omega
_{1}-\left( p_{112}+\frac{1}{3}\left( 3s+2\right) Kp_{1}+\frac{1}{3}%
sK_{1}u\right) \omega _{2}, \\
\varepsilon _{12} &=&dp_{12}-\left( s+2\right) p_{12}~\gamma -\left( p_{112}-%
\frac{1}{6}\left( 3s+2\right) Kp_{1}-\frac{1}{6}sK_{1}u\right) \omega _{1} \\
&&-\left( p_{122}+\frac{1}{6}\left( 3s+2\right) Kp_{2}+\frac{1}{6}%
sK_{2}u\right) \omega _{2}, \\
\varepsilon _{22} &=&dp_{22}-\left( s+2\right) p_{22}~\gamma -\left( p_{122}-%
\frac{1}{3}\left( 3s+2\right) Kp_{2}-\frac{1}{3}sK_{2}u\right) \omega
_{1}-p_{222}\omega _{2}.
\end{eqnarray*}%
In a similar way, from the relations%
\begin{eqnarray*}
\delta _{1}\delta _{112} &=&\delta _{1112}-\frac{1}{6}\left( 3s+4\right)
K\delta _{11}-\frac{1}{6}\left( 2s+1\right) K_{1}\delta _{1}-\frac{1}{12}%
sK_{11}, \\
\delta _{1}\delta _{122} &=&\delta _{1122}-\frac{1}{3}\left( 3s+4\right)
K\delta _{12}-\frac{1}{6}\left( 2s+1\right) K_{2}\delta _{1}-\frac{1}{6}%
\left( 2s+1\right) K_{1}\delta _{2}-\frac{1}{6}sK_{12}, \\
\delta _{1}\delta _{222} &=&\delta _{1222}-\frac{1}{2}\left( 3s+4\right)
K\delta _{22}-\frac{1}{2}\left( 2s+1\right) K_{2}\delta _{2}-\frac{1}{4}%
sK_{22}, \\
\delta _{2}\delta _{111} &=&\delta _{1112}+\frac{1}{2}\left( 3s+4\right)
K\delta _{11}+\frac{1}{2}\left( 2s+1\right) K_{1}\delta _{1}+\frac{1}{4}%
sK_{11}, \\
\delta _{2}\delta _{112} &=&\delta _{1122}+\frac{1}{3}\left( 3s+4\right)
K\delta _{12}+\frac{1}{6}\left( 2s+1\right) K_{2}\delta _{1}+\frac{1}{6}%
\left( 2s+1\right) K_{1}\delta _{2}+\frac{1}{6}sK_{12}, \\
\delta _{2}\delta _{122} &=&\delta _{1222}+\frac{1}{6}\left( 3s+4\right)
K\delta _{22}+\frac{1}{6}\left( 2s+1\right) K_{2}\delta _{2}+\frac{1}{12}%
sK_{22},
\end{eqnarray*}%
we get the following representation for the fourth-order Cartan forms:%
\begin{eqnarray*}
&&\varepsilon _{111}=dp_{111}-\left( s+3\right) p_{111}\gamma
-p_{1111}\omega _{1} \\
&&-\left( p_{1112}+\frac{1}{2}\left( 3s+4\right) Kp_{11}+\frac{1}{2}\left(
2s+1\right) K_{1}p_{1}+\frac{1}{4}sK_{11}u\right) \omega _{2}, \\
&&\varepsilon _{112}=dp_{112}-\left( s+3\right) p_{112}\gamma \\
&&-\left( p_{1112}-\frac{1}{6}\left( 3s+4\right) Kp_{11}-\frac{1}{6}\left(
2s+1\right) K_{1}p_{1}-\frac{1}{12}sK_{11}u\right) \omega _{1} \\
&&-\left( p_{1122}+\frac{1}{3}\left( 3s+4\right) Kp_{12}+\frac{1}{6}\left(
2s+1\right) K_{2}p_{1}+\frac{1}{6}\left( 2s+1\right) K_{1}p_{2}+\frac{1}{6}%
sK_{12}u\right) \omega _{2}, \\
&&\varepsilon _{122}=dp_{122}-\left( s+3\right) p_{122}\gamma \\
&&-\left( p_{1122}-\frac{1}{3}\left( 3s+4\right) Kp_{12}-\frac{1}{6}\left(
2s+1\right) K_{2}p_{1}-\frac{1}{6}\left( 2s+1\right) K_{1}p_{2}-\frac{1}{6}%
sK_{12}u\right) \omega _{1} \\
&&-\left( p_{1222}+\frac{1}{6}\left( 3s+4\right) Kp_{22}+\frac{1}{6}\left(
2s+1\right) K_{2}p_{2}+\frac{1}{12}sK_{22}u\right) \omega _{2}, \\
&&\varepsilon _{222}=dp_{222}-\left( s+3\right) p_{222}\gamma \\
&&-\left( p_{1222}-\frac{1}{2}\left( 3s+4\right) Kp_{22}-\frac{1}{2}\left(
2s+1\right) K_{2}p_{2}-\frac{1}{4}sK_{22}u\right) \omega _{1}-p_{2222}\omega
_{2}.
\end{eqnarray*}

\subsection{The Total Derivative and the Mayer Bracket}

We shall denote by $\widehat{X}$ the total derivative corresponding to a
vector field $X$ on the manifold $M^{2}$ (see, for example, \cite{KLV} or %
\cite{ALV}). Using the representations of Cartan's forms, we get the
following expressions for the vector fields $\widehat{\partial _{1}}$ and $%
\widehat{\partial _{2}}$:%
\begin{eqnarray*}
&&\widehat{\partial _{1}}=\partial _{1}+\left( sg_{1}u+p_{1}\right) \frac{%
\partial }{\partial u}+\left( \left( s+1\right) g_{1}p_{1}+p_{11}\right) 
\frac{\partial }{\partial p_{1}}+\left( \left( s+1\right) g_{1}p_{2}+p_{12}-%
\frac{s}{2}Ku\right) \frac{\partial }{\partial p_{2}} \\
&&+\left( \left( s+2\right) g_{1}p_{11}+p_{111}\right) \frac{\partial }{%
\partial p_{11}}+\left( \left( s+2\right) g_{1}p_{12}+p_{112}-\frac{3s+4}{6}%
Kp_{1}-\frac{s}{6}K_{1}u\right) \frac{\partial }{\partial p_{12}} \\
&&+\left( \left( s+2\right) g_{1}p_{22}+p_{122}-\frac{3s+4}{3}Kp_{2}-\frac{s%
}{3}K_{2}u\right) \frac{\partial }{\partial p_{22}} \\
&&+\left( \left( s+3\right) g_{1}p_{111}+p_{1111}\right) \frac{\partial }{%
\partial p_{111}} \\
&&+\left( \left( s+3\right) g_{1}p_{112}+p_{1112}-\frac{3s+4}{6}Kp_{11}-%
\frac{2s+1}{6}K_{1}p_{1}-\frac{s}{12}K_{11}u\right) \frac{\partial }{%
\partial p_{112}} \\
&&+\left( \left( s+3\right) g_{1}p_{122}+p_{1122}-\frac{3s+4}{3}Kp_{12}-%
\frac{2s+1}{6}K_{2}p_{1}-\frac{2s+1}{6}K_{1}p_{2}-\frac{s}{6}K_{12}u\right) 
\frac{\partial }{\partial p_{122}} \\
&&+\left( \left( s+3\right) g_{1}p_{222}+p_{1222}-\frac{3s+4}{2}Kp_{22}-%
\frac{2s+1}{2}K_{2}p_{2}-\frac{s}{4}K_{22}u\right) \frac{\partial }{\partial
p_{222}}+\cdots
\end{eqnarray*}%
and%
\begin{eqnarray*}
&&\widehat{\partial _{2}}=\partial _{2}+\left( sg_{2}u+p_{2}\right) \frac{%
\partial }{\partial u}+\left( \left( s+1\right) g_{2}p_{1}+p_{12}+\frac{s}{2}%
Ku\right) \frac{\partial }{\partial p_{1}}+\left( \left( s+1\right)
g_{2}p_{2}+p_{22}\right) \frac{\partial }{\partial p_{2}} \\
&&+\left( \left( s+2\right) g_{2}p_{11}+p_{112}+\frac{3s+4}{3}Kp_{1}+\frac{s%
}{3}K_{1}u\right) \frac{\partial }{\partial p_{11}} \\
&&+\left( \left( s+2\right) g_{2}p_{12}+p_{122}+\frac{3s+4}{6}Kp_{2}+\frac{s%
}{6}K_{2}u\right) \frac{\partial }{\partial p_{12}} \\
&&+\left( \left( s+2\right) g_{2}p_{22}+p_{222}\right) \frac{\partial }{%
\partial p_{22}} \\
&&+\left( \left( s+3\right) g_{2}p_{111}+p_{1112}+\frac{3s+4}{2}Kp_{11}+%
\frac{2s+1}{2}K_{1}p_{1}+\frac{s}{4}K_{11}u\right) \frac{\partial }{\partial
p_{111}} \\
&&+\left( \left( s+3\right) g_{2}p_{112}+p_{1122}+\frac{3s+4}{3}Kp_{12}+%
\frac{2s+1}{6}K_{2}p_{1}+\frac{2s+1}{6}K_{1}p_{2}+\frac{s}{6}K_{12}u\right) 
\frac{\partial }{\partial p_{112}} \\
&&+\left( \left( s+3\right) g_{2}p_{122}+p_{1222}+\frac{3s+4}{6}Kp_{22}+%
\frac{2s+1}{6}K_{2}p_{2}+\frac{s}{12}K_{22}u\right) \frac{\partial }{%
\partial p_{122}} \\
&&+\left( \left( s+3\right) g_{2}p_{222}+p_{2222}\right) \frac{\partial }{%
\partial p_{222}}+\cdots
\end{eqnarray*}%
\newline

We shall denote by 
\begin{equation*}
\widehat{\delta _{i}}\left( h\right) =\widehat{\partial _{i}}(h)-wg_{i}h
\end{equation*}%
the covariant derivatives of a function $h$ of weight $w$ on the jet space $%
\mathbb{J}^{r}\left( s\right) $ and call it the \emph{total covariant
derivative }of $h$ \emph{\ }along $\partial _{i}.$ As earlier, we shall
denote the symmetrized total derivatives by $\widehat{\delta }%
_{i_{1}...i_{l}}.$

In these notations, the linearization of a function $h$ of weight $w$ on $%
\mathbb{J}^{r}\left( s\right) $ (cf. \cite{KLV} and \cite{ALV}) has the form%
\begin{equation*}
l_{h}=\sum_{(i_{1}...i_{l})}\frac{\partial ^{l}h}{\partial p_{i_{1}...i_{l}}}%
~\widehat{\delta }_{i_{1}...i_{l}},
\end{equation*}%
and the Mayer bracket (see \cite{KL}) of functions $f$ and $g$ of weights $%
w_{1}$ and $w_{2}$ defined correspondingly on $\mathbb{J}^{n}\left( s\right) 
$ and $\mathbb{J}^{m}\left( s\right) $ has the form 
\begin{equation*}
\lbrack f,g]=\sum_{(i_{1}...i_{n})}\frac{\partial ^{n}f}{\partial
p_{i_{1}...i_{n}}}~\widehat{\delta }_{i_{1}...i_{n}}\left( g\right)
-\sum_{(j_{1},...,j_{m})}\frac{\partial ^{m}f}{\partial p_{j_{1}...j_{m}}}~%
\widehat{\delta }_{j_{1}...j_{m}}\left( g\right) .
\end{equation*}

\section{The Mayer Bracket and the First Obstruction for Linearizability}

Let us rewrite equations (\ref{mu-invariants}) symbolically. The functions
on $\mathbb{J}^{2}\left( 1\right) ,$ that correspond to these equations, are 
\begin{eqnarray*}
I_{1} &=&p_{11}-2p_{12}-up_{1}+2up_{2}+K_{1}, \\
I_{2} &=&p_{22}-2p_{12}-2up_{1}+up_{2}+K_{2}.
\end{eqnarray*}%
Equations (\ref{AG-compatibility eq}) are compatible if and only if the
Mayer bracket of $I_{1}$ and $I_{2}\ $ vanishes (see \cite{KL}). In our case,%
\begin{equation*}
I_{12}=[I_{1},I_{2}]=\widehat{\delta }_{11}\left( I_{2}\right) -\widehat{%
\delta }_{22}\left( I_{1}\right) +2\widehat{\delta }_{12}\left(
I_{1}-I_{2}\right)
\end{equation*}%
or%
\begin{eqnarray*}
&&I_{12}=u\left( -2p_{111}+3p_{112}+3p_{122}-2p_{222}\right) \\
&&+3\left( p_{2}-2p_{1}\right) p_{11}+6\left( p_{1}+p_{2}\right)
p_{12}+3\left( p_{1}-2p_{2}\right) p_{22}+8K(p_{11}-p_{12}+p_{22}) \\
&&+3(2K_{1}-K_{2})p_{1}+3(2K_{2}-K_{1})p_{2}+u\left(
K_{11}-K_{12}+K_{22}\right) \\
&&+3\left( K_{112}-K_{122}\right) .\ \ 
\end{eqnarray*}%
Solving the first prolongation of the system%
\begin{eqnarray*}
\widehat{\delta }_{1}\left( I_{1}\right) &=&0,\ \ \widehat{\delta }%
_{2}\left( I_{1}\right) =0, \\
\widehat{\delta }_{1}\left( I_{2}\right) &=&0,\ \ \ \widehat{\delta }%
_{2}\left( I_{2}\right) =0
\end{eqnarray*}%
with respect to $p_{111},\;p_{112,}\;p_{122}$ and $p_{222}$ and substituting
the result into $I_{12},$ we get%
\begin{eqnarray*}
I_{12} &=&24Kp_{12}+6(2K_{1}-K_{2})p_{1}+6(2K_{2}-K_{1})p_{2}+24Ku\left(
p_{1}-p_{2}\right) \\
&&+3u\left( K_{11}-K_{12}+K_{22}\right) -8K\left( K_{1}+K_{2}\right)
+3\left( K_{112}-K_{122}\right) -3Ku^{3}.
\end{eqnarray*}

Note that%
\begin{equation*}
\tau \left( I_{12}\right) =I_{12}.
\end{equation*}%
Solving the equations 
\begin{equation*}
I_{1}=0,\;I_{2}=0,\;I_{12}=0
\end{equation*}%
with respect to $p_{ij},$ we obtain 
\begin{equation}  \label{PijSystem}
\renewcommand{\arraystretch}{1.3} 
\begin{array}{lll}
12Kp_{11} & = & 3Ku^{3}+12Kup_{1}+6(K_{2}-2K_{1})p_{1}+6(K_{1}-2K_{2})p_{2}
\\ 
&  & -3u\left( K_{11}-K_{12}+K_{22}\right) +3\left( K_{122}-K_{112}\right)
+4K(2K_{2}-K_{1}), \\ 
12Kp_{22} & = & 3Ku^{3}+12Kup_{2}+6(K_{2}-2K_{1})p_{1}+6(K_{1}-2K_{2})p_{2}
\\ 
&  & -3u\left( K_{11}-K_{12}+K_{22}\right) +3\left( K_{122}-K_{112}\right)
+4K(2K_{1}-K_{2}), \\ 
24Kp_{12} & = & 
3Ku^{3}+24Ku(p_{2}-p_{1})+6(K_{2}-2K_{1})p_{1}+6(K_{1}-2K_{2})p_{2} \\ 
&  & -3u\left( K_{11}-K_{12}+K_{22}\right) +3\left( K_{122}-K_{112}\right)
+8K(K_{1}+K_{2}).%
\end{array}
\renewcommand{\arraystretch}{1}
\end{equation}%
The expressions for the symmetric covariant derivatives of the curvature
function \thinspace $K$ are given in section 10.1.

We write down the above equations in the form 
\begin{eqnarray*}
p_{11} &=&P_{11}\left( u,p_{1},p_{2},K\right) , \\
p_{12} &=&P_{12}\left( u,p_{1},p_{2},K\right) , \\
p_{22} &=&P_{22}\left( u,p_{1},p_{2},K\right) .
\end{eqnarray*}%
In order to find their compatibility, first, taking $s=1,$ we derive from (%
\ref{delta 1-2 products})$~$ that 
\begin{eqnarray*}
\delta _{2}\delta _{11}-\delta _{1}\delta _{12} &=&\frac{5}{2}K\delta _{1}+%
\frac{1}{2}K_{1}, \\
\delta _{2}\delta _{12}-\delta _{1}\delta _{22} &=&\frac{5}{2}K\delta _{2}+%
\frac{1}{2}K_{2}.
\end{eqnarray*}%
It follows that the equations 
\begin{eqnarray*}
\widehat{\delta }_{2}\left( P_{11}\right) -\widehat{\delta }_{1}\left(
P_{12}\right) -\frac{5}{2}Kp_{1}-\frac{1}{2}K_{1}u &=&0, \\
\widehat{\delta }_{2}\left( P_{12}\right) -\widehat{\delta }_{1}\left(
P_{22}\right) -\frac{5}{2}Kp_{2}-\frac{1}{2}K_{2}u &=&0
\end{eqnarray*}%
are the compatibility conditions for (\ref{PijSystem}).

Let us denote by $G_{1}$ and $G_{2}$ the left-hand sides of the above
equations into which the values of $p_{ij}$ taken from (\ref{PijSystem}) are
substituted.

These functions are polynomials in $p_{i}$ and $u$ of the form%
\begin{eqnarray*}
G_{1} &=&(p_{1}^{2}-2p_{1}p_{2})+A_{11}p_{1}+A_{12}p_{2}+A_{10}, \\
G_{2} &=&(p_{2}^{2}-2p_{1}p_{2})+A_{21}p_{1}+A_{22}p_{2}+A_{20},
\end{eqnarray*}%
where all coefficients are functions of the curvature $K$ and its covariant
derivatives up to order four:%
\begin{eqnarray*}
A_{11} &=&\frac{5u^{2}}{8}+\frac{3u\left( K_{1}-K_{2}\right) }{4K}-\frac{13K%
}{4}-\frac{7(K_{1}-2K_{2})(2K_{1}-K_{2})}{16K^{2}} \\
&&+\frac{7(5K_{11}+5K_{22}-11K_{12})}{8K},\  \\
\ A_{12} &=&-\frac{5u^{2}}{4}-\frac{3u\ K_{1}}{4K}+\frac{7(K_{1}-2K_{2})^{2}%
}{16K^{2}}+\frac{5K_{12}-2K_{11}-5K_{22}}{4K}, \\
A_{10} &=&\frac{9u^{3}(K_{1}-2K_{2})}{96K}+u\left( -\frac{5K_{1}}{2}+\frac{%
21(2K_{2}-K_{1})(K_{11}-K_{12}+K_{22})}{96K^{2}}\right) \\
&&+\frac{u(K_{111}-2K_{222}-3K_{112}+3K_{122})}{8K}+\frac{%
9(2K_{22}-2K_{12}-K_{11})}{16} \\
&&-\frac{21(2K_{2}-K_{1})(K_{122}-K_{112})}{96K^{2}}+\frac{%
17(2K_{2}^{2}-2K_{1}K_{2}-K_{1}^{2})}{48K} \\
&&+\frac{2K_{1222}-3K_{1122}+K_{1112}}{8K},
\end{eqnarray*}%
and 
\begin{equation*}
A_{21}=-\tau \left( A_{12}\right) ,A_{22}=-\tau \left( A_{11}\right)
,A_{20}=\tau \left( A_{10}\right)
\end{equation*}%
because%
\begin{equation*}
\tau \left( G_{1}\right) =G_{2}.
\end{equation*}

The following theorem outlines the successive steps in the investigation of
solvability for main equations (\ref{AG2-equations}).

\begin{theorem}
\label{1-st Golds Th} \ 

\begin{enumerate}
\item Differential equations $(\ref{AG2-equations})$ are solvable with
respect to the functions $\lambda _{1}$ and $\lambda _{2}$ if and only if
the function $\mu $ satisfies differential equations $(\ref{AG-compatibility
eq})$.

\item For the system of differential equations $(\ref{AG-compatibility eq})$
be solvable, one needs to add the compatibility condition $I_{12}=0$ to this
system.

\item The compatibility conditions for the resulting system $(\ref{PijSystem}%
)$ have the form 
\begin{equation}
G_{1}=0,\ G_{2}=0.  \label{G12 system}
\end{equation}
\end{enumerate}
\end{theorem}

\section{The Second Obstruction for Linearizability}

In this section, we investigate the solvability of the system of equations (%
\ref{PijSystem}) and (\ref{G12 system}). To this end, we differentiate the
left-hand sides of (\ref{G12 system}), 
\begin{eqnarray*}
G_{11} &=&\widehat{\delta _{1}}\left( G_{1}\right) ,\ G_{12}^{s}=\frac{1}{2}%
\left( \widehat{\delta _{1}}\left( G_{2}\right) +\widehat{\delta _{2}}\left(
G_{1}\right) \right) , \\
\ G_{22} &=&\widehat{\delta _{2}}\left( G_{2}\right) ,\,G_{12}^{a}=\frac{1}{2%
}\left( \widehat{\delta _{1}}\left( G_{2}\right) -\widehat{\delta _{2}}%
\left( G_{1}\right) \right)
\end{eqnarray*}%
and substitute the second covariant derivatives taken from (\ref{PijSystem})
into the result of differentiation.

Finally, we arrive at the system 
\begin{equation}
G_{1}=0,\,G_{2}=0,\;G_{11}=0,\ G_{12}^{s}=0,\;G_{12}^{a}=0,\ G_{22}=0,
\label{GijSystem}
\end{equation}%
which is equivalent to system (\ref{PijSystem})--(\ref{G12 system}).

By the construction, we get the symmetry 
\begin{equation*}
\tau \left( G_{11}\right) =G_{22},\ \tau \left( G_{12}^{s}\right)
=G_{12}^{s},\ \tau \left( G_{12}^{a}\right) =-G_{12}^{a},\ \tau \left(
G_{22}\right) =G_{11}.
\end{equation*}%
In the coordinates, these functions have the form%
\begin{equation}  \label{Gijexpressions}
\renewcommand{\arraystretch}{2.0} 
\begin{array}{lll}
G_{11} & = & -\displaystyle\frac{\left( K_{1}+K_{2}\right) }{4K}p_{1}^{2}+%
\frac{7K_{1}-8K_{2}}{4K}p_{1}p_{2}+\frac{2K_{2}-K_{1}}{K}p_{2}^{2} \\ 
&  & +A_{111}~p_{1}+A_{112}~p_{2}+A_{110}+\displaystyle\frac{5u}{4}G_{1}, \\ 
G_{12}^{s} & = & \displaystyle \frac{8K_{1}-7K_{2}}{8K}p_{1}^{2}+\frac{%
K_{1}+K_{2}}{2K}p_{1}p_{2}-\frac{7K_{1}-8K_{2}}{8K}p_{2}^{2} \\ 
&  & +A_{121}p_{1}+A_{122}p_{2}+A_{120}+\displaystyle\frac{5u}{4}G_{1}-\frac{%
5u}{4}G_{2}, \\ 
G_{12}^{a} & = & \displaystyle\frac{39}{4}%
up_{1}p_{2}+B_{121}p_{1}+B_{122}p_{2}+B_{120} \\ 
&  & +\left( \displaystyle\frac{13}{4}u-\frac{3K_{2}}{8K}\right)
G_{1}+\left( \frac{13}{4}u+\frac{3K_{1}}{8K}\right) G_{2}, \\ 
G_{22} & = & \displaystyle\frac{2K_{1}-K_{2}}{K}p_{1}^{2}+\frac{7K_{2}-8K_{1}%
}{4K}p_{1}p_{2}-\frac{K_{1}+K_{2}}{4K}p_{2}^{2} \\ 
&  & +A_{221}~p_{1}+A_{222}~p_{2}+A_{220}-\displaystyle\frac{5u}{4}G_{2},%
\end{array}
\renewcommand{\arraystretch}{1}
\end{equation}%
where%
\begin{eqnarray*}
A_{111} &=&\frac{3}{32}u^{3}-\frac{3\left( 7K_{1}-12K_{2}\right) }{32K}%
u^{2}+\cdots , \\
A_{112} &=&-\frac{3}{16}u^{3}+\frac{3K_{1}}{16}u^{2}+\cdots , \\
A_{110} &=&-\frac{3\left( K_{1}-2K_{2}\right) }{128K}u^{4}+\frac{%
33K_{1}\left( 2K_{2}-K_{1}\right) }{128K^{2}}u^{3}+\frac{3\left(
2K_{12}-K_{11}\right) }{16K}u^{3}+\cdots ,
\end{eqnarray*}%
and%
\begin{eqnarray*}
A_{121} &=&-\frac{3}{32}u^{3}+\frac{3\left( 5K_{2}-14K_{1}\right) }{64K}%
u^{2}+\cdots , \\
A_{122} &=&-\frac{3}{64}u^{3}+\frac{3\left( 14K_{2}-5K_{1}\right) }{128K}%
u^{2}+\cdots , \\
A_{120} &=&\frac{3\left( K_{1}+K_{2}\right) }{128K}u^{4}+\frac{33\left(
K_{2}^{2}-K_{1}^{2}\right) }{128K^{2}}u^{3}+\frac{3\left(
K_{22}-K_{11}\right) }{16K}u^{2}+\cdots
\end{eqnarray*}%
and%
\begin{eqnarray*}
B_{121} &=&-\frac{195}{16}u^{3}+\frac{9\left( 9K_{2}-5K_{1}\right) }{8K}%
u^{2}+\cdots , \\
B_{122} &=&-\frac{195}{16}u^{3}+\frac{9\left( 9K_{1}-5K_{2}\right) }{8K}%
u^{2}+\cdots , \\
B_{120} &=&\frac{15}{32}u^{5}+\frac{117\left( K_{2}-K_{1}\right) }{64K}%
u^{4}+\cdots
\end{eqnarray*}%
Moreover, 
\begin{equation*}
A_{222}=-\tau \left( A_{111}\right) ,\;\,A_{221}=-\tau \left( A_{112}\right)
,\;\,A_{220}=\tau \left( A_{110}\right) .
\end{equation*}

The detailed expressions for these coefficients can be found in Section 10.2.

Summarizing, we get the following system of first-order PDEs on the function 
$\mu $: 
\begin{equation}
G_{1}=0,\;G_{2}=0,\;G_{11}=0,\ G_{12}^{s}=0,\ G_{12}^{a}=0,\,G_{22}=0,
\label{Gi&GijSystem}
\end{equation}%
which is equivalent to system (\ref{AG-compatibility eq}).

We remark that this system is symmetric with respect to the involution $\tau
.$

Next we note that equations (\ref{Gi&GijSystem}) contain only linear
combinations of the functions $p_{1},\,p_{2},\,p_{1}^{2},\,p_{1}p_{2},%
\,p_{2}^{2}$ with coefficients depending on $u,\,K$ and the covariant
derivatives of $K\ $up to order five.

We solve the equations $G_{1}=0,\ G_{12}^{a}=0,\ G_{2}=0$ with respect to $%
p_{1}^{2},\;p_{1}p_{2},\;p_{2}^{2}$.

The determinant of the system is equal to $39u/4$.

\ Note that $\mu =0$ implies $K_{1}=K_{2}=0$ due to (\ref{AG-compatibility
eq}), and it is impossible for nonparallelizable $3$-webs.

Indeed, if $K_{1}=K_{2}=0,$ then $\partial _{1}\left( K\right) =\partial
_{1}\left( K\right) =2HK,$ and 
\begin{equation*}
0=H\left( \partial _{2}-\partial _{1}\right) \left( K\right) =[\partial
_{1},\partial _{2}]\left( K\right) =2K\left( \partial _{1}\left( H\right)
-\partial _{2}\left( H\right) \right) =-2K^{2}.
\end{equation*}

Solving the equations $G_{1}=0,\ G_{12}^{a}=0,\ G_{2}=0$ with respect to $%
p_{1}^{2},\,p_{1}p_{2},\,p_{2}^{2},$ we get the expressions for $p_{i}p_{j}$
in the form of linear combinations of $p_{1}$ and $p_{2}$.

Substituting these expressions into the system $G_{11}=0,\;G_{22}=0$ and
solving the resulting system of linear equations with respect to $p_{1}$ and 
$p_{2},$ we find that%
\begin{equation*}
p_{1}=\frac{V_{1}}{V_{0}},\ p_{2}=\frac{V_{2}}{V_{0}},
\end{equation*}%
where $V_{1}$ and $V_{2}$ are polynomials of degree eight with respect to $%
u, $ and their coefficients depend on the curvature function $K$ and its
covariant derivatives up to order five. The leading terms of $V_{1}$ and $%
V_{2}$ are%
\begin{eqnarray*}
V_{1} &=&-\frac{3^{4}}{2^{8}}KK_{1}u^{8}+\frac{3^{2}}{2^{5}}%
[7(K_{1}^{2}+2K_{1}K_{2}-2K_{2}^{2})+13K(-K_{11}-2K_{12}+2K_{22})]u^{7}+%
\cdots , \\
V_{2} &=&-\frac{3^{4}}{2^{8}}KK_{2}u^{8}+\frac{3^{2}}{2^{5}}%
[7(2K_{1}^{2}-2K_{1}K_{2}-K_{2}^{2})+13K(K_{22}+2K_{12}-2K_{11})]u^{7}+%
\cdots ,
\end{eqnarray*}%
and the denominator $V_{0}$ is the seven-degree polynomial (see section 10.3)%
\begin{equation*}
V_{0}=-\frac{13\cdot 3^{3}}{2^{6}}K^{2}u^{7}+\frac{3^{3}}{2^{5}}%
[15(K_{1}^{2}-K_{1}K_{2}+K_{2}^{2})+13K(K_{11}-K_{12}+K_{22})]u^{5}+\cdots
\end{equation*}%
As we have seen, the functions $p_{i}p_{j}\;$are linear combinations of $%
p_{1}$ and $p_{2}$. Substituting the above expressions for $p_{1}$ and $%
p_{2} $ into the expressions for $p_{i}p_{j},$ we get%
\begin{equation*}
p_{1}^{2}=\frac{V_{11}}{V_{0}},\;p_{1}p_{2}=\frac{V_{12}}{V_{0}},\ p_{2}^{2}=%
\frac{V_{22}}{V_{0}},
\end{equation*}%
where $V_{ij}$ are polynomials of degree $11$ with respect to $u$ and their
coefficients depend on the curvature function $K$ and its covariant
derivatives up to order five. The leading terms of $V_{ij}\;$are%
\begin{eqnarray*}
V_{11} &=&\frac{5\cdot 3^{3}}{2^{9}}K^{2}u^{11}-\frac{3^{6}}{2^{10}}%
KK_{1}u^{10} \\
&&+\frac{3^{2}}{2^{10}}%
[35K_{1}^{2}+412K_{1}K_{2}-412K_{2}^{2}+20K(-16K_{11}-23K_{12}+23K_{22})]u^{9}+\cdots ,
\\
V_{12} &=&\frac{5\cdot 3^{3}}{2^{10}}K^{2}u^{11}+\frac{3^{6}}{2^{10}}K\left(
K_{2}-K_{1}\right) u^{10} \\
&&+\frac{3^{2}}{2^{10}}%
[206(K_{1}^{2}-K_{2}^{2})+653K_{1}K_{2}+10K(23K_{11}-101K_{12}+23K_{22})]u^{9}+\cdots ,
\\
V_{22} &=&\frac{5\cdot 3^{3}}{2^{9}}K^{2}u^{11}+\frac{3^{6}}{2^{10}}%
KK_{2}u^{10} \\
&&-\frac{3^{2}}{2^{10}}%
[412K_{1}^{2}-412K_{1}K_{2}-35K_{2}^{2}+20K(-23K_{11}+23K_{12}+16K_{22})]u^{9}+\cdots
\end{eqnarray*}

Note that the equation $G_{12}^{s}=0$ holds automatically.

The resulting system%
\begin{equation}  \label{pipj^2System}
\renewcommand{\arraystretch}{2.0} 
\begin{array}{lll}
p_{1} & = & \displaystyle\frac{V_{1}}{V_{0}},\ \ p_{2}=\frac{V_{2}}{V_{0}},
\\ 
p_{1}^{2} & = & \displaystyle\frac{V_{11}}{V_{0}},\ p_{1}p_{2}=\frac{V_{12}}{%
V_{0}},\ \ p_{2}^{2}=\frac{V_{22}}{V_{0}}%
\end{array}
\renewcommand{\arraystretch}{1}
\end{equation}%
is $\tau $-symmetric:%
\begin{eqnarray*}
\tau \left( V_{0}\right) &=&-V_{0},\ \tau \left( V_{1}\right) =V_{2},\;\tau
\left( V_{2}\right) =V_{1},\  \\
\tau \left( V_{11}\right) &=&-V_{22},\;\tau \left( V_{12}\right)
=-V_{12},\;\tau \left( V_{22}\right) =-V_{11}.
\end{eqnarray*}

This system gives us the following polynomial equations on $u$:%
\begin{equation}
V_{0}~V_{11}-V_{1}^{2}=0,\ V_{0}~V_{22}-V_{2}^{2}=0,\
V_{0}~V_{12}-V_{1}V_{2}=0.  \label{PolynomEquations}
\end{equation}%
Let us denote the left-hand sides of the above equations by $Q_{ij}$ and $%
Q_{a}\;$and $Q_{s}$ symmetrizations of $Q_{11}$\ and $Q_{22}.$ We consider
the polynomials 
\begin{eqnarray*}
2Q_{a} &=&Q_{11}+Q_{22}=V_{0}~(V_{11}-V_{22})-V_{1}^{2}+V_{2}^{2}, \\
2Q_{s} &=&Q_{11}-Q_{22}=V_{0}~(V_{11}+V_{22})-V_{1}^{2}-V_{2}^{2} \\
\ Q_{12} &=&V_{0}~V_{12}-V_{1}V_{2}.
\end{eqnarray*}%
The degree of each of the polynomials $Q_{s}$\ and $Q_{12}$\ equals $18$
while the degree of $Q_{a}\;$does not exceed $17$:

\begin{eqnarray*}
Q_{a} &=&\frac{13\cdot 3^{9}}{2^{17}}K^{3}\left( K_{1}+K_{2}\right) ~u^{17}-%
\frac{3^{6}K^{2}}{2^{16}}[973\left( K_{1}^{2}-K_{2}^{2}\right) +1690K\left(
K_{22}-K_{11}\right) ]~u^{16}+\cdots , \\
Q_{s} &=&-\frac{65\cdot 3^{6}}{2^{15}}K^{4}u^{18}+\frac{13\cdot 3^{9}}{2^{17}%
}K^{3}\left( K_{1}-K_{2}\right) ~u^{17} \\
&&-\frac{3^{5}}{2^{16}}K^{2}[-3337(K_{1}^{2}+K_{2}^{2})+6256K_{1}K_{2}+130K%
\left( K_{11}-40K_{12}+K_{22}\right) ]~u^{16}+\cdots , \\
Q_{12} &=&-\frac{65\cdot 3^{6}}{2^{16}}K^{4}u^{18}+\frac{13\cdot 3^{9}}{%
2^{16}}K^{3}(K_{1}-K_{2})~u^{17} \\
&&-\frac{243}{2^{15}}%
[K^{2}(-1564(K_{1}^{2}+K_{2}^{2})+4483K_{1}K_{2}+65K(10K_{11}-49K_{12}+10K_{22})]~u^{16}+\cdots ,
\end{eqnarray*}

In order to complete integration of system (\ref{pipj^2System}), we
differentiate one of equations (\ref{PolynomEquations}), say, the first one, 
\begin{equation}  \label{TotalDerPolynom}
\renewcommand{\arraystretch}{2.2} 
\begin{array}{lll}
\displaystyle\frac{\partial Q_{a}}{\partial u}p_{1}+\widehat{\delta _{1}^{K}}%
\left( Q_{a}\right) & = & 0, \\ 
\displaystyle\frac{\partial Q_{a}}{\partial u}p_{2}+\widehat{\delta _{2}^{K}}%
\left( Q_{a}\right) & = & 0,%
\end{array}
\renewcommand{\arraystretch}{1}
\end{equation}%
where $\widehat{\delta _{i}^{K}}$ are the total derivatives relative to $K$
(see 10.4 for the expressions of $\widehat{\delta _{i}^{K}}$):%
\begin{eqnarray*}
\widehat{\delta _{1}^{K}} &=&K_{1}\frac{\partial }{\partial K}+K_{11}\frac{%
\partial }{\partial K_{1}}+\left( K_{12}-K^{2}\right) \frac{\partial }{%
\partial K_{2}}+\cdots \\
\widehat{\delta _{2}^{K}} &=&K_{2}\frac{\partial }{\partial K}+(K_{12}+K^{2})%
\frac{\partial }{\partial K_{1}}+K_{22}\frac{\partial }{\partial K_{2}}%
+\cdots
\end{eqnarray*}

Substituting the covariant derivatives $p_{1}$ and $p_{2}$ taken from the
first two equations of (\ref{pipj^2System}) into (\ref{TotalDerPolynom}), we
get the new system of polynomial equations on $u$: 
\begin{eqnarray*}
\frac{\partial Q_{a}}{\partial u}V_{1}+V_{0}\ \widehat{\delta _{1}^{K}}%
\left( Q_{a}\right) &=&0, \\
\frac{\partial Q_{a}}{\partial u}V_{2}+V_{0}\ \widehat{\delta _{2}^{K}}%
\left( Q_{a}\right) &=&0.
\end{eqnarray*}

The polynomials\vspace{1pt}%
\begin{equation*}
Q_{1}=\frac{\partial Q_{a}}{\partial u}V_{1}+V_{0}\widehat{~\delta _{1}^{K}}%
\left( Q_{a}\right)
\end{equation*}%
and%
\begin{equation*}
Q_{2}=\frac{\partial Q_{a}}{\partial u}V_{2}+V_{0}\ \widehat{\delta _{2}^{K}}%
\left( Q_{a}\right)
\end{equation*}%
are of degree $24,$ and their coefficients depend on the curvature function $%
K$ and its covariant derivatives up to order six:%
\begin{eqnarray*}
Q_{1} &=&\frac{131\cdot 65\cdot 3^{9}}{2^{23}}K^{5}K_{1}u^{24}+\cdots , \\
Q_{2} &=&\frac{131\cdot 65\cdot 3^{9}}{2^{23}}K^{5}K_{2}u^{24}+\cdots
\end{eqnarray*}

The next result follows from the above consideration and is basic for
finding linearizability conditions for $3$-webs.

\begin{theorem}
Let $W$ be a nonparallelizable $3$-web. Then the smooth solvability of the
system of nonlinear partial differential equations 
\begin{equation*}
I_{1}\left( \mu \right) =0,\ I_{2}\left( \mu \right) =0
\end{equation*}%
is equivalent to the existence of real and smooth solutions of the following
system of algebraic equations$:$%
\begin{equation*}
Q_{a}=0,\ Q_{s}=0,\ Q_{12}=0,\ Q_{1}=0,\ Q_{2}=0.
\end{equation*}
\end{theorem}

In 1912 Gronwall (\cite{gr})\ made the following conjecture: \textit{if a
nonparallelizable }$3$\textit{-web }$W_{3}$\textit{\ in the plane is
linearizable, then, up to a projective transformation, a diffeomorphism
transforming }$W_{3}$\textit{\ into a linear }$3$\textit{-web, is uniquely
determined. }

Bol (\cite{bol1},\cite{bol2}, 1938) and Bor\r{u}vka (\cite{bo}, 1938) proved
that the number of projectively nonequivalent linearizations of a
nonparallelizable linearizable $3$-web does not exceed $16.$ Vaona (\cite{v2}%
, 1961) reduced this number to $11.$ Grifone,\ Muzsnay and Saab (\cite{GMS
01}, 2001) proved that this number does not exceed $15.$

The above theorem implies the following result.

\begin{corollary}
Let $W$ be a nonparallelizable, linearizable $3$-web. Then the number of
projectively nonequivalent linearizations of such a web does not exceed $15.$
\end{corollary}

\begin{proof}
Observe, that if $\mu $ satisfies the system $I_{1}\left( \mu \right) =0,\
I_{2}\left( \mu \right) =0,$ then system (\ref{AG2-equations}) is completely
integrable, and its solutions $\left( \lambda _{1},\lambda _{2}\right) $ are
determined by values $\lambda _{1}\left( a_{0}\right) $ and $\lambda
_{2}\left( a_{0}\right) $ at some fixed point $a_{0}\in M.$ Moreover, it is
easy to check that the projective transformations act transitively on the
set of $\left( \lambda _{1}\left( a_{0}\right) ,\lambda _{2}\left(
a_{0}\right) \right) .$ So, up to a projective transformation, the values $%
\left( \lambda _{1}\left( a_{0}\right) ,\lambda _{2}\left( a_{0}\right)
\right) $ are nonessential.

As we showed earlier, the polynomials $Q_{a},\;Q_{s}$ and $Q_{12}$ are of
degrees $17,\;18,\;$and$\;18,$ and each of the polynomials $Q_{1}\;$and$%
\;Q_{2}$ is of degree $24.$ Hence, there is a linear combination $L\;$of $%
Q_{s}$ and $Q_{12}$ having degree$\ \leq 17,$ and there is a linear
combination $S\;$of $Q_{a}$ and $L\;$having degree $\leq 16.$

In fact, we can take as $L$ the \ polynomial%
\begin{eqnarray*}
L &=&Q_{s}-2Q_{12}=\frac{13\cdot 3^{10}}{2^{17}}K^{3}(K_{2}-K_{1})u^{17} \\
&&+\frac{3^{6}}{2^{15}}%
K^{2}[973(4K_{1}K_{2}-K_{1}^{2}-K_{2}^{2})+1690K(K_{11}-4K_{12}+K_{22})]u^{16}+\cdots
\end{eqnarray*}%
If $K_{2}-K_{1}\neq 0\;$and $K_{1}+K_{2}\neq 0,$\ then as\ $S$ we can take
the polynomial%
\begin{eqnarray*}
S &=&(K_{1}+K_{2})L-3(K_{2}-K_{1})Q_{a} \\
&=&\frac{3^{6}}{2^{15}}%
K^{2}[-1946(K_{1}^{3}+K_{2}^{3})+2919K_{1}K_{2}(K_{1}+K_{2})+1690K(2K_{1}-K_{2})K_{11}
\\
&&-3380K(K_{1}+K_{2})K_{12}+1690K(2K_{2}-K_{1})K_{22}]u^{16}+\cdots
\end{eqnarray*}%
If $K_{2}-K_{1}=0\;$(or $K_{1}+K_{2}=0$),$\;$then the polynomial $L\;$(resp. 
$Q_{a}$) is already of degree $16.$

Thus the polynomials $Q_{a},\;Q_{s},$ $Q_{12}$ and $Q_{1},\;Q_{2}$ can have
at most $16$ common roots. One of these roots gives $\mu $ for the $3$-web
under consideration. Therefore, the number of projectively nonequivalent
linearizations of the web $W$ does not exceed $15$.
\end{proof}

\textbf{Remark.} In the paper \cite{AGL} we have proved \ that $\mu $ is
uniquely determined by the basic invariant of linearizable $d$-webs, if $%
d\geq 4.$ The above proof shows that the Gronwall conjecture is correct for
such webs. Namely, up to a projective transformation, for linearizable $d$%
-webs, $d\geq 4,$ there exists a unique linearization.

\section{Differential Invariants for Linearizability and the Blaschke
Conjecture}

In this section we consider the case of nonparallelizable, linearizable $3$%
-webs. We will need some new algebraic constructions.

\subsection{Resultant and Its Generalizations}

Let $T,S_{1},...,S_{n}$ be polynomials over an algebraically closed field $%
\mathbb{F},$ $T,S_{1},...,S_{n}\in \mathbb{F[}u],$ and $\limfunc{char}%
\mathbb{F}=0.$ Denote by $\mathbf{R}(f,g)$ the resultant of polynomials $f$
and $g.$ Recall that $\mathbf{R}(f,g)$ as a function in $g$ given $f$ is
homogeneous of degree $\deg f.$ Hence $\mathbf{R}(T,x_{1}S_{1}+x_{2}S_{2}+%
\cdots +x_{n}S_{n})$ is a homogeneous polynomial of degree $\deg T$ in $%
x_{1},...,x\,_{n}$: 
\begin{equation*}
\mathbf{R}(T,\sum_{i=1}^{n}x_{i}S_{i})=\sum_{\sigma }x^{\sigma }\mathbf{R}%
_{\sigma }(T,S_{1},...,S_{n}),
\end{equation*}%
where $\sigma $ runs over all multi-indices of the length $\deg T,$i.e.,

\begin{equation*}
\ \mathbf{R}(T,%
\sum_{i=1}^{n}x_{i}S_{i})=x_{1}^{i_{1}}x_{2}^{i_{2}}...x_{n}^{i_{n}}\mathbf{R%
}_{i_{1}i_{2}...i_{n}}(T,S_{1},...,S_{n}).
\end{equation*}%
We call the coefficients $\mathbf{R}_{\sigma }(T,S_{1},...,S_{n})$ $($\emph{%
generalized}$\emph{)}$\emph{\ resultants }of the system of polynomials $%
T,S_{1},...,S_{n}.$

\begin{theorem}
The polynomials $T,S_{1},...,S_{n}$ have a common root if and only if all
resultants $\mathbf{R}_{\sigma }(T,S_{1},...,S_{n})$ are equal to zero.
\end{theorem}

\begin{proof}
To illustrate the idea of the proof and to avoid unnecessary technicalities,
we consider only the case $n=2.$ Assume also that the leading coefficient of 
$T$ is equal to $1.$

Let $\lambda _{1},..,\lambda _{t}$ be roots of $T,$ and $t=\deg T.$ Then 
\begin{equation*}
\mathbf{R}(T,x_{1}S_{1}+x_{2}S_{2})=\dprod\limits_{i=1}^{t}\left(
x_{1}S_{1}\left( \lambda _{i}\right) +x_{2}S_{2}\left( \lambda _{i}\right)
\right) =\sum_{a=0}^{t}x_{1}^{a}x_{2}^{t-a}\mathbf{R}_{a,t-a}(T,S_{1},S_{2}),
\end{equation*}%
where $\mathbf{R}_{t,0}(T,S_{1},S_{2})=\mathbf{R}(T,S_{1}),\mathbf{R}%
_{0,t}(T,S_{1},S_{2})=\mathbf{R}(T,S_{2}),$ and for $1\leq a\leq t-1,$ we
get 
\begin{equation*}
\mathbf{R}_{a,t-a}(T,S_{1},S_{2})=\sum_{I}S_{1}\left( \lambda
_{i_{1}}\right) \cdots S_{1}\left( \lambda _{i_{a}}\right) S_{2}\left(
\lambda _{j_{1}}\right) \cdots S_{2}\left( \lambda _{j_{t-a}}\right) .
\end{equation*}%
Here we have denoted by $\left( j_{1},...,j_{t-a}\right) $ the multi-index
complementary to $I=\left( i_{1},...,i_{a}\right) .$

First, let $T,S_{1}\;$and$\;S_{2}$ have a common root. Then the polynomials $%
T$ and $x_{1}S_{1}+x_{2}S_{2}$ have a common root for all $x_{1},x_{2},$ and
therefore $\mathbf{R}_{a,t-a}(T,S_{1},S_{2})=0$ for all $a.$ Conversely, let 
$\mathbf{R}_{a,t-a}(T,S_{1},S_{2})=0$ for all $a.$ Then $\mathbf{R}(T,S_{1})=%
\mathbf{R}_{t,0}(T,S_{1},S_{2})=0,\ \ \mathbf{R}(T,S_{2})=\mathbf{R}%
_{0,t}(T,S_{1},S_{2})=0,$ and therefore $T$ and $S_{1}$ have a common root,
say $\nu ,$ and $T$ and $S_{2}$ have a common root, say $\mu .$ Assume that
they have no more common roots, and consider, for example, $\mathbf{R}%
_{1,t-1}(T,S_{1},S_{2}).$

One has 
\begin{equation*}
\mathbf{R}_{1,t-1}(T,S_{1},S_{2})=S_{1}\left( \mu \right) \cdot S_{2}\left(
\nu \right) \cdot S_{2}\left( \lambda _{j_{1}}\right) \cdot ...\cdot
S_{2}\left( \lambda _{j_{t-2}}\right) =0,
\end{equation*}%
where $\left( \lambda _{1},..,\lambda _{t}\right) =\nu \cup \mu \cup \left(
\lambda _{1},..,\lambda _{t-2}\right) $ is the disjoint union.

Hence, either $S_{1}\left( \mu \right) =0$ or $S_{2}\left( \nu \right) =0,$
and therefore $T,S_{1}$ and $S_{2}$ have a common root.

In the case when the polynomials have\ common roots of multiplicity two or
higher, $\mathbf{R}_{1,t-1}(T,S_{1},S_{2})=\mathbf{R}%
_{1-t,1}(T,S_{1},S_{2})=0,$ and vanishing of $\mathbf{R}%
_{2,t-2}(T,S_{1},S_{2})$ shows that $T,S_{1}$ and $S_{2}$ have a common
root, etc.
\end{proof}

\textbf{Remark. }The number of resultants $\mathbf{R}_{\sigma
}(T,S_{1},...,S_{n})$ equals the dimension of homogeneous polynomials of
degree $t=\deg T$ in $n$\ variables, and therefore equals 
\begin{equation*}
\binom{n+t-1}{t}.
\end{equation*}

\subsection{Differential Invariants for Linearizability}

As we have seen earlier, the solvability of the system of differential
equations (\ref{AG-compatibility eq}) is equivalent to the existence of real
roots of the system of algebraic equations%
\begin{equation}
Q_{a}=0,\ Q_{s}=0,\ Q_{12}=0,\ Q_{1}=0,\ Q_{2}=0.  \label{ALGEBRAIC SYSTEM}
\end{equation}%
We apply the above theorem and get the following result.

\begin{theorem}
Let $W$ be a nonparallelizable $3$-web. If the $3$-web $W$ is linearizable,
then the following differential invariants%
\begin{equation*}
\mathbf{R}_{i_{1}i_{2}i_{3}i_{4}}\left(
Q_{a},Q_{s},Q_{12},Q_{1},Q_{2}\right) 
\end{equation*}%
vanish, and algebraic system $(\ref{ALGEBRAIC SYSTEM})$ has at least one
real smooth solution.

Conversely, if the differential invariants vanish and algebraic system $(\ref%
{ALGEBRAIC SYSTEM})$ has at least one real smooth solution, then the $3$-web
is linearizable.
\end{theorem}

Note that all the differential invariants depend on the curvature function $K
$ and its covariant derivatives up to order six, but $\mathbf{R}%
_{i_{1}i_{2}i_{3}i_{4}}\left( Q_{a},Q_{s},Q_{12},Q_{1},Q_{2}\right) $ \ with 
$i_{3}=i_{4}=0$ depend on the curvature function $K$ and its covariant
derivatives up to order five. Since for a nonparallelizable $3$-web, we have 
$\deg \;Q_{a}=17,$ the total number of invariants equals $1040=\binom{4+17-1%
}{17},$ and among them there are $18=\binom{2+17-1}{17}$ invariants of order
five in $K.$ In terms of the web function $f\left( x,y\right) ,$ the
corresponding orders are nine and eight.

Note also that the number of invariants is not invariant: it depends which
of the polynomials $Q_{a},Q_{s},Q_{12},Q_{1},Q_{2}\;$we take as the first
one. In our considerations we took the polynomial $Q_{a}$ of the least
degree $17\;$as the first polynomial. Moreover, the number of invariants can
be reduced if we find a linear combination of the above five polynomials
whose degree is less than $17$, replace one of five polynomials by this
linear combination and take this combination as the first polynomial (see
our earlier considerations where we found a polynomial of degree not
exceeding $16$).

\textbf{Remark. }In the book \cite{B 55}$\;(\S 17)$ Blaschke made the
following conjecture: The linearizability conditions for a nonparallelizable 
$3$-web are expressed in terms of the web function$\;f\left( x,y\right) \;$%
and its covariant derivatives up to order nine, and the table in $\S 17$
shows that the number of differential invariants equals four. As we have
seen, Blaschke's estimate of the ''functional codimension'' of the orbits of
the linearizable $3$-webs was correct while the number of algebraic
conditions \ is much greater than four. Moreover, not all linearizability
invariants are of order nine: eighteen of them are of order eight.

To find out whether algebraic system (\ref{ALGEBRAIC SYSTEM}) has real
solutions, we consider the greatest common divisor $\mathbf{G}=\mathbf{GCD}%
[Q_{a},Q_{s},Q_{12},Q_{1},Q_{2}]$ of the polynomials $%
Q_{a},Q_{s},Q_{12},Q_{1},Q_{2}.$

The following theorem, which is important when one is testing a $3$-web for
linearizability, is obvious.

\begin{theorem}
If $\deg \mathbf{G=}0,$ then there are no common solutions, and the $3$-web
is nonlinearizable. If $\deg \mathbf{G>}1,$ but $\mathbf{G}$ has no real
roots, then the $3$-web is also nonlinearizable. In the case when $\deg 
\mathbf{G=}1,$ or $\deg \mathbf{G>}$\ $1$ but $\mathbf{G}$ has a real root,
a $3$-web is linearizable.
\end{theorem}

Note that in the latter case, the number of real roots can give us an
improvement of our estimate of the Gronwall conjecture: if the number of
real roots of $\mathbf{G}$ equals $s,$\ then \ the number of projectively
nonequivalent linearizations of a nonparallelizable, nonexceptional
linearizable $3$-web$\;W$ does not exceed\ $s$. If $s<15,$\thinspace then
this will be an improvement of our estimate.

\subsection{Linear 3-Webs}

To test our systems of equations and the differential invariants, we
consider them for linear $3$-webs.

Let us assume that the web function $f(x,y)$ defines a linear $3$-web, and
let $d_{\text{st}}$ be the covariant differential of the flat connection in
coordinates $x$ and $y.$

Then $d_{\text{st}}\left( a~dx\right) =dx\otimes da$ and $d_{\text{st}%
}\left( a~dy\right) =dy\otimes da.$ Therefore,%
\begin{equation*}
d_{\text{st}}\left( \omega _{1}\right) =\omega _{1}\otimes \frac{df_{x}}{%
f_{x}},\ \ d_{\text{st}}\left( \omega _{2}\right) =\omega _{2}\otimes \frac{%
df_{y}}{f_{y}},
\end{equation*}%
and the affine deformation tensor $T=d_{\gamma }-d_{\text{st}}$ between the
Chern and the flat connections equals to%
\begin{eqnarray*}
T(\omega _{1}) &=&\left( \frac{f_{xx}}{f_{x}^{2}}-H\right) \omega
_{1}\otimes \omega _{1}, \\
T(\omega _{2}) &=&\left( \frac{f_{yy}}{f_{y}^{2}}-H\right) \omega
_{2}\otimes \omega _{2}.
\end{eqnarray*}%
Therefore, for linear $3$-webs we have 
\begin{eqnarray*}
\lambda _{1} &=&\lambda _{2}=0,\; \\
\mu &=&\frac{f_{xx}}{f_{x}^{2}}-H=-\left( \frac{f_{yy}}{f_{y}^{2}}-H\right) .
\end{eqnarray*}

If we assume that $f\left( x,y\right) $ is a solution of the Euler equation,
i.e., $f_{x}=f~f_{y}$, then we get 
\begin{equation*}
\mu =\frac{1}{f}.
\end{equation*}

Moreover, in this case%
\begin{equation*}
H=\frac{1}{f}+\frac{f_{yy}}{f_{y}^{2}},
\end{equation*}%
and the curvature function is%
\begin{equation*}
K=-\frac{f_{yy}}{f~f_{y}^{2}}.
\end{equation*}%
The first covariant derivatives of the curvature function are%
\begin{eqnarray*}
K_{1} &=&-\frac{2K}{f}+\frac{f_{yyy}}{ff_{y}^{3}}, \\
K_{2} &=&-\frac{K}{f}+\frac{f_{yyy}}{ff_{y}^{3}}
\end{eqnarray*}%
and%
\begin{equation*}
K_{2}-K_{1}=\frac{K}{f}.
\end{equation*}%
Note that the covariant derivatives of the function $\mu $ are 
\begin{equation*}
\delta _{1}\left( \mu \right) =\delta _{2}\left( \mu \right) =K.
\end{equation*}%
One can check that equations (\ref{AG2-equations}) and (\ref%
{AG-compatibility eq}) have the common solution $\mu =1/f,$ and the same is
true for equations (\ref{Gi&GijSystem}).

Moreover $1/f$ is a common real root for algebraic system (\ref{ALGEBRAIC
SYSTEM}).

\subsection{Procedure for Applying the Linearizability Criterion\ \ \ \ \ \
\ \ \ \ \ \ \ \ \ \ \ \ \ \ \ \ \ \ \ \ \ \ \ \ \ \ \ \ \ \ \ \ \ \ \ \ \ \
\ \ \ \ \ \ \ \ \ \ \ \ \ }

Now we can outline a procedure which can be applied to determine whether a $%
3 $-web $W_{3}\;$given by a web function $z=f(x,y)\;$is linearizable:

\begin{enumerate}
\item Compute the curvature $K$ and its covariant derivatives up to order
five (see formula (\ref{curvature main}) and formulas in Section 10.1).

\item Compute $A_{ij},\;A_{ijk,}$ and $B_{ijk,}\;i,j,k=0,1,2\;$(see formulas
in Sections 8 and 10.2).

\item Compute the polynomial $V_{0}\;$(see Sections 8 and 10.3).

\item Compute the polynomials $V_{ij},\;V_{i}$ and $%
Q_{a},Q_{s},Q_{12},Q_{1},Q_{2}\;$(see Sections 8 and 10.3).

\item Compute $\mathbf{G}=\mathbf{GCD}[Q_{a},Q_{s},Q_{12},Q_{1},Q_{2}]$ and
apply the linearizability condition outlined in Theorem 9.3.
\end{enumerate}

\subsection{Examples}

\textbf{Example 6} We consider the $3$-web in the plane with the web function%
\begin{equation*}
f(x,y)=x^{2}+xy+y^{2}.
\end{equation*}

For this web we have:%
\begin{eqnarray*}
H &=&\frac{1}{(2x+y)(x+2y)},\;K=-\frac{6(x^{2}-y^{2})}{(2x+y)^{3}(x+2y)^{3}},
\\
K_{1} &=&-\frac{6(4x^{3}+3x^{2}y-12xy^{2}-13y^{3})}{(2x+y)^{5}(x+2y)^{4}}, \\
\;K_{2} &=&-\frac{6(13x^{3}+12x^{2}y-3xy^{2}-4y^{3})}{(2x+y)^{4}(x+2y)^{5}}%
,...,
\end{eqnarray*}%
and%
\begin{eqnarray*}
Q_{a} &=&\frac{13\cdot 3^{14}(x^{2}-y^{2})^{4}(5x^{2}-8xy+5y^{2})}{%
2^{12}(2x+y)^{14}(x+2y)^{14}}\mu ^{17}+..., \\
Q_{s} &=&-\frac{65\cdot 3^{10}(x^{2}-y^{2})^{4}}{2^{11}(2x+y)^{12}(x+2y)^{12}%
}\mu ^{18}+..., \\
Q_{12} &=&-\frac{65\cdot 3^{10}(x^{2}-y^{2})^{4}}{%
2^{12}(2x+y)^{12}(x+2y)^{12}}\mu ^{18}+..., \\
Q_{1} &=&\frac{13\cdot 3^{14}(x^{2}-y^{2})^{5}}{2^{12}(2x+y)^{22}(x+2y)^{21}}%
\cdot (5700x^{5}+13577x^{4}y-2480x^{3}y^{2} \\
&&-37710x^{2}y^{3}-44660xy^{4}-18343y^{5})\mu ^{24}+..., \\
Q_{2} &=&\frac{13\cdot 3^{14}(x^{2}-y^{2})^{5}}{2^{12}(2x+y)^{21}(x+2y)^{22}}%
\cdot (18343x^{5}+44660x^{4}y+37710x^{3}y^{2} \\
&&+2480x^{2}y^{3}-13577xy^{4}-5700y^{5})\mu ^{24}+...
\end{eqnarray*}%
Evaluating the polynomials $Q_{a},\;Q_{s},\;Q_{12}\;$and $Q_{1},\;Q_{2}\;$at
the point ($0.1,1),\;$we find that%
\begin{eqnarray*}
F_{a} &=&Q_{a}(0.1,1)=0.204819\mu ^{17}+..., \\
F_{s} &=&Q_{s}(0.1,1)=-0.0274492\mu ^{18}+..., \\
F_{12} &=&Q_{12}(0.1,1)=-0.0137246\mu ^{18}+..., \\
F_{1} &=&Q_{1}(0.1,1)=-3.94038\mu ^{24}+..., \\
F_{2} &=&Q_{1}(0.1,1)=-0.678834\mu ^{24}+...
\end{eqnarray*}%
We calculate now the resultant of the polynomials $F_{a}\;$and$\;F_{12}$:%
\begin{equation*}
\mathbf{R}(F_{a},F_{12})=-1.046\cdot 10^{185}\neq 0.
\end{equation*}%
Since the resultant of $F_{a}\;$and$\;F_{12}\;$does not vanish, the
polynomials $Q_{a}\;$and$\;Q_{12}\;$(and therefore the polynomials $%
Q_{a},Q_{s},Q_{12},Q_{1},Q_{2}$)$\;$have no common roots, and as a result, 
\emph{the }$3$\emph{-web under consideration is not linearizable.}

\textbf{Remark. }Note that even if the resultants of all pairs of the
polynomials $Q_{a},Q_{s},Q_{12},Q_{1},Q_{2}$\ were vanished, we could not
make any conclusion---the further investigation involving the generalized
resultants or finding the greatest common divisor $\mathbf{G}=\mathbf{GCD}%
[Q_{a},Q_{s},Q_{12},Q_{1},Q_{2}]$ would be necessary to answer the question
whether the $3$-web under consideration is linearizable or not linearizable.

\textbf{Example 7} We consider the $3$-web in the plane with the web function%
\begin{equation*}
f(x,y)=(x+y)e^{-x}
\end{equation*}

$($see Example $5)$.\ 

For this web we have:%
\begin{eqnarray*}
H &=&\frac{e^{x}}{x+y-1},\ \ K=\frac{e^{2x}}{(1-x-y)^{3}},\;\  \\
K_{1} &=&\frac{3e^{3x}}{(x+y-1)^{5}},\;\;K_{2}=-\frac{e^{3x}}{(x+y-1)^{4}}%
,...,
\end{eqnarray*}%
and%
\begin{eqnarray*}
Q_{a} &=&\frac{13\cdot 3^{9}e^{9x}(x+y-4)}{2^{17}(x+y-1)^{12}}\mu ^{17}+...,
\\
Q_{s} &=&-\frac{65\cdot 3^{6}e^{8x}}{2^{15}(x+y-1)^{12}}\mu ^{18}+..., \\
Q_{12} &=&-\frac{65\cdot 3^{6}e^{8x}}{2^{16}(x+y-1)^{12}}\mu ^{18}+..., \\
Q_{1} &=&\frac{13\cdot 3^{12}e^{14x}(829x+829y-3472)}{41\cdot 31\cdot
11\cdot 3\cdot 2^{3}(x+y-1)^{22}}\mu ^{24}+..., \\
Q_{2} &=&-\frac{13\cdot 3^{12}e^{14x}(259x+259y-1192)}{41\cdot 31\cdot
11\cdot 3\cdot 2^{3}(x+y-1)^{21}}\mu ^{24}+...
\end{eqnarray*}%
Evaluating the polynomials $Q_{a},\;Q_{s},\;Q_{12}\;$and $Q_{1},\;Q_{2}\;$at
the point ($0,0.1),\;$we find that%
\begin{eqnarray*}
F_{a} &=&Q_{a}(0,0.1)=-33.2808\mu ^{17}+..., \\
F_{s} &=&Q_{s}(0,0.1)=-5.12013\mu ^{18}+..., \\
F_{12} &=&Q_{12}(0,0.1)=-2.56006\mu ^{18}+..., \\
F_{1} &=&Q_{1}(0,0.1)=-7085.94\mu ^{24}+..., \\
F_{2} &=&Q_{1}(0,0.1)=-2194.28\mu ^{24}+...
\end{eqnarray*}%
We calculate now the resultant of the polynomials $F_{11}\;$and$\;F_{22}$:%
\begin{equation*}
\mathbf{R}(F_{a},F_{12})=1.23007\cdot 10^{272}\neq 0.
\end{equation*}%
Since the resultant of $F_{a}\;$and$\;F_{12}\;$does not vanish, the
polynomials $Q_{a}\;$and$\;Q_{12}\;$(and therefore the polynomials $%
Q_{a},Q_{s},Q_{12},Q_{1},Q_{2}$)$\;$have no common roots, and as a result, 
\emph{the }$3$\emph{-web under consideration is not linearizable.}

\section{Appendix. Computational Formulae}

\subsection{Symmetrized Covariant Derivatives of the Curvature}

{\ {\ 
\begin{equation*}
K_{12}=\delta _{1}\left( K_{2}\right) +K^{2}.
\end{equation*}
}}

\begin{eqnarray*}
K_{112} &=&\delta _{1}\left( K_{12}\right) +\frac{5}{3}KK_{1}; \\
K_{122} &=&\delta _{1}\left( K_{22}\right) +\frac{10}{3}KK_{2}.
\end{eqnarray*}

\begin{eqnarray*}
K_{1112} &=&\delta _{1}\left( K_{112}\right) +\frac{11}{6}KK_{11}+\frac{5}{6}%
K_{1}^{2}; \\
K_{1122} &=&\delta _{1}\left( K_{122}\right) +\frac{11}{3}KK_{12}+\frac{5}{3}%
K_{1}K_{2}; \\
K_{1222} &=&\delta _{1}\left( K_{222}\right) +\frac{11}{2}KK_{22}+\frac{5}{2}%
K_{2}^{2}.
\end{eqnarray*}

\begin{eqnarray*}
K_{11112} &=&\delta _{1}\left( K_{1112}\right) +\frac{21}{10}\left(
KK_{111}+K_{1}K_{11}\right) ; \\
K_{11122} &=&\delta _{1}\left( K_{1122}\right) +\frac{7}{5}\left(
3KK_{112}+K_{2}K_{11}+2K_{1}K_{12}\right) ; \\
K_{11222} &=&\delta _{1}\left( K_{1222}\right) +\frac{21}{10}\left(
3KK_{122}+K_{1}K_{22}+2K_{2}K_{12}\right) ; \\
K_{12222} &=&\delta _{1}\left( K_{2222}\right) +\frac{42}{5}\left(
KK_{222}+K_{2}K_{22}\right) .
\end{eqnarray*}

\begin{eqnarray*}
K_{111112} &=&\delta _{1}\left( K_{11112}\right) +\frac{12}{5}KK_{1111}+%
\frac{14}{5}K_{1}K_{111}+\frac{7}{5}K_{11}^{2}; \\
K_{111122} &=&\delta _{1}\left( K_{11122}\right) +\frac{1}{5}\left(
24KK_{1112}+7K_{2}K_{111}+21K_{1}K_{112}+14K_{12}K_{11}\right) ; \\
K_{111222} &=&\delta _{1}\left( K_{11222}\right) \\
&&+\frac{1}{5}\left(
36KK_{1122}+21K_{1}K_{122}+21K_{2}K_{112}+7K_{11}K_{22}+14K_{12}^{2}\right) ;
\\
K_{112222} &=&\delta _{1}\left( K_{12222}\right) +\frac{2}{5}\left(
24KK_{1222}+7K_{1}K_{222}+21K_{2}K_{122}+14K_{12}K_{22}\right) ; \\
K_{122222} &=&\delta _{1}\left( K_{22222}\right)
+12KK_{2222}+14K_{2}K_{222}+7K_{22}^{2}.
\end{eqnarray*}%
\newline

\subsection{Coefficients $\emph{A}_{ijk}$ and $\emph{B}_{ijk}\;$in ( \ref%
{Gijexpressions})}

\begingroup\scalefont{0.8}

Here we give the expressions of the coefficients $A_{ijk}\;$and $B_{ijk}\;$%
in formulas (\ref{Gijexpressions}) for $G_{11},\;G_{12}^{s},\;G_{12}^{a},\;$%
and$\;G_{22}\;$(see Section 8):

\begin{eqnarray*}
&&A_{111}= \\
&&\frac{3u^{3}}{32}+\frac{3\left( -7K_{1}+12K_{2}\right) u^{2}}{32K}+\frac{%
145Ku}{16}+\frac{\left( 26K_{1}^{2}-95K_{1}K_{2}-10K_{2}^{2}\right) u}{%
64K^{2}} \\
&&-\frac{\left( 13K_{11}-43K_{12}+13K_{22}\right) u}{32K}+\frac{1}{48}\left(
-179K_{1}-62K_{2}\right) \\
&&+\frac{77K_{1}\left( 2K_{1}^{2}-5K_{1}K_{2}+2K_{2}^{2}\right) }{64K^{3}}+%
\frac{K_{2}\left( 45K_{11}-41K_{12}+7K_{22}\right) }{16K^{2}} \\
&&-\frac{K_{1}\left( 95K_{11}-145K_{12}+27K_{22}\right) }{32K^{2}}+\frac{%
3K_{111}-8K_{112}+5K_{122}-K_{222}}{4K};
\end{eqnarray*}%
\begin{eqnarray*}
&&A_{112}= \\
&&-\frac{3u^{3}}{16}+\frac{3K_{1}u^{2}}{16K}+\frac{\left(
53K_{1}^{2}-20K_{1}K_{2}+20K_{2}^{2}\right) u}{64K^{2}}-\frac{\left(
2K_{11}+13K_{12}-13K_{22}\right) u}{16K} \\
&&-\frac{31}{24}\left( K_{1}-2K_{2}\right) -\frac{77K_{1}\left(
K_{1}-2K_{2}\right) ^{2}}{64K^{3}}-\frac{17K_{2}\left( K_{11}-2K_{12}\right) 
}{8K^{2}} \\
&&+\frac{K_{1}\left( 25K_{11}-54K_{12}+20K_{22}\right) }{16K^{2}}-\frac{%
2K_{111}-7K_{112}+7K_{122}}{4K};
\end{eqnarray*}%
\begin{eqnarray*}
&&A_{110}= \\
&&\frac{3\left( K_{1}-2K_{2}\right) u^{4}}{128K}\text{{}}-\frac{%
33K_{1}\left( K_{1}-2K_{2}\right) u^{3}}{128K^{2}}+\frac{3\left(
K_{11}-2K_{12}\right) u^{3}}{16K} \\
&&-\frac{23K_{2}\left( K_{11}-K_{12}+K_{22}\right) u^{2}}{64K^{2}}+\frac{%
23K_{1}\left( 16K^{2}+K_{11}-K_{12}+K_{22}\right) u^{2}}{128K^{2}} \\
&&-\frac{5\left( K_{111}-3K_{112}+3K_{122}-2K_{222}\right) u^{2}}{32K}-\frac{%
31}{64}\left( K_{11}+2K_{12}-2K_{22}\right) u \\
&&+\frac{77K_{1}\left( K_{1}-2K_{2}\right) \left(
K_{11}-K_{12}+K_{22}\right) u}{128K^{3}}-\frac{5\left( K_{11}-2K_{12}\right)
\left( K_{11}-K_{12}+K_{22}\right) u}{16K^{2}} \\
&&+\frac{K_{2}\left( 28K_{111}-51K_{112}+51K_{122}\right) u}{64K^{2}}-\frac{%
21K_{1}^{2}-198K_{1}K_{2}+198K_{2}^{2}}{64K}u \\
&&+\frac{K_{1}\left( -44K_{111}+99K_{112}-99K_{122}+32K_{222}\right) u}{%
128K^{2}} \\
&&+\frac{8K_{1111}-34K_{1112}+54K_{1122}-36K_{1222}}{64K}u+\text{{}}\frac{11%
}{3}K\left( K_{1}-2K_{2}\right) + \\
&&\frac{17K_{1}\left( 11K_{1}^{2}-20K_{1}K_{2}+20K_{2}^{2}\right) }{192K^{2}}%
-\frac{353K_{2}\left( K_{11}-2K_{12}\right) }{480K} \\
&&+\frac{K_{1}\left( -907K_{11}+398K_{12}-1104K_{22}\right) }{960K}-\frac{%
5\left( K_{11}-2K_{12}\right) \left( K_{112}-K_{122}\right) }{16K^{2}} \\
&&+\frac{77K_{1}\left( K_{1}-2K_{2}\right) \left( K_{112}-K_{122}\right) }{%
128K^{3}}-\frac{33}{40}\left( K_{111}-K_{112}+K_{122}\right) \\
&&+\frac{7K_{2}\left( K_{1112}-K_{1122}\right) }{16K^{2}}+\frac{K_{1}\left(
-11K_{1112}+19K_{1122}-8K_{1222}\right) }{32K^{2}} \\
&&+\frac{K_{11112}-3K_{11122}+2K_{11222}}{8K};
\end{eqnarray*}

\begin{eqnarray*}
&&A_{221}= \\
&&-\frac{3u^{3}}{16}-\frac{3K_{2}u^{2}}{16K}+\frac{%
20K_{1}^{2}-20K_{1}K_{2}+53K_{2}^{2}}{64K^{2}}u+\frac{%
13K_{11}-13K_{12}-2K_{22}}{16K}u \\
&&-\frac{31}{24}\left( -2K_{1}+K_{2}\right) +\frac{77K_{2}\left(
-2K_{1}+K_{2}\right) ^{2}}{64K^{3}}+\frac{17K_{1}\left(
-2K_{12}+K_{22}\right) }{8K^{2}} \\
&&-\frac{K_{2}\left( 20K_{11}-54K_{12}+25K_{22}\right) }{16K^{2}}+\frac{%
7K_{112}-7K_{122}+2K_{222}}{4K};
\end{eqnarray*}

\begin{eqnarray*}
&&A_{222}= \\
&&\frac{3u^{3}}{32}\text{{}}+\frac{3\left( -12K_{1}+7K_{2}\right) }{32K}%
u^{2}-\frac{10K_{1}^{2}+95K_{1}K_{2}-26K_{2}^{2}}{64K^{2}}u \\
&&-\frac{145Ku}{16}-\frac{13K_{11}-43K_{12}+13K_{22}}{32K}u+\frac{1}{48}%
\left( -62K_{1}-179K_{2}\right) \\
&&-\frac{77K_{2}\left( 2K_{1}^{2}-5K_{1}K_{2}+2K_{2}^{2}\right) }{64K^{3}}-%
\frac{K_{1}\left( 7K_{11}-41K_{12}+45K_{22}\right) }{16K^{2}} \\
&&+\frac{K_{2}\left( 27K_{11}-145K_{12}+95K_{22}\right) }{32K^{2}}+\frac{%
K_{111}-5K_{112}+8K_{122}-3K_{222}}{4K};
\end{eqnarray*}

\begin{eqnarray*}
&&A_{220}= \\
&&-\frac{3\left( -2K_{1}+K_{2}\right) u^{4}}{128K}+\frac{33K_{2}\left(
-2K_{1}+K_{2}\right) u^{3}}{128K^{2}}+\frac{6K_{12}-3K_{22}}{16K}u^{3} \\
&&-\frac{23K_{2}u^{2}}{8}+\frac{23\left( -2K_{1}+K_{2}\right) \left(
K_{11}-K_{12}+K_{22}\right) u^{2}}{128K^{2}} \\
&&-\frac{5\left( -2K_{111}+3K_{112}-3K_{122}+K_{222}\right) u^{2}}{32K}+%
\frac{31}{64}\left( 2K_{11}-2K_{12}-K_{22}\right) u \\
&&-\frac{3\left( 66K_{1}^{2}-66K_{1}K_{2}+7K_{2}^{2}\right) u}{64K}-\frac{%
5\left( 2K_{12}-K_{22}\right) \left( K_{11}-K_{12}+K_{22}\right) u}{16K^{2}}
\\
&&-\frac{77K_{2}\left( 2K_{1}\left( K_{12}-K_{22}\right) +K_{2}\left(
K_{11}+K_{22}\right) \right) u}{128K^{3}}+\frac{%
7^{2}11^{2}K_{1}K_{2}^{3}K_{11}K_{12}u}{2^{13}K^{6}} \\
&&-\frac{K_{1}\left( 51K_{112}-51K_{122}+28K_{222}\right) u}{64K^{2}}+\frac{%
18K_{1112}-27K_{1122}+17K_{1222}-4K_{2222}}{32K}u \\
&&+\frac{K_{2}\left( -32K_{111}+99K_{112}-99K_{122}+44K_{222}\right) u}{%
128K^{2}}+\frac{11}{3}K\left( -2K_{1}+K_{2}\right) \\
&&-\frac{17K_{2}\left( 20K_{1}^{2}-20K_{1}K_{2}+11K_{2}^{2}\right) }{192K^{2}%
}+\frac{353K_{1}\left( -2K_{12}+K_{22}\right) }{480K} \\
&&+\frac{K_{2}\left( 1104K_{11}-398K_{12}+907K_{22}\right) }{960K}+\frac{%
77\left( 2K_{1}-K_{2}\right) K_{2}\left( K_{112}-K_{122}\right) }{128K^{3}}
\\
&&-\frac{5\left( 2K_{12}K_{112}-2K_{12}K_{122}+K_{22}K_{122}\right) }{16K^{2}%
}+\frac{33\left( K_{112}-K_{122}+K_{222}\right) }{40} \\
&&+\frac{5K_{22}K_{112}}{16K^{2}}-\frac{K_{2}\left(
8K_{1112}-19K_{1122}+11K_{1222}\right) }{32K^{2}}+\frac{14K_{1}\left(
-K_{1122}+K_{1222}\right) }{32K^{2}} \\
&&+\frac{2K_{11122}-3K_{11222}+K_{12222}}{8K};
\end{eqnarray*}

\begin{eqnarray*}
&&A_{121}=-\frac{3u^{3}}{32}+\frac{3\left( -14K_{1}+5K_{2}\right) }{64K}%
u^{2}+\frac{145Ku}{16}+\frac{10K_{1}^{2}-79K_{1}K_{2}+61K_{2}^{2}}{64K^{2}}u
\\
&&+\frac{13K_{11}+17K_{12}-17K_{22}}{32K}u+\frac{77\left(
2K_{1}^{3}-K_{1}^{2}K_{2}-2K_{1}K_{2}^{2}+K_{2}^{3}\right) }{64K^{3}} \\
&&+\frac{1}{96}\left( 124K_{1}-365K_{2}\right) +\frac{K_{2}\left(
55K_{11}+79K_{12}-81K_{22}\right) }{64K^{2}} \\
&&+\frac{K_{1}\left( -190K_{11}+54K_{12}+64K_{22}\right) }{64K^{2}}+\frac{%
3K_{111}-K_{112}-2K_{122}+K_{222}}{4K};
\end{eqnarray*}

\begin{eqnarray*}
A_{122} &=&-\frac{3u^{3}}{32}+\frac{3\left( -5K_{1}+14K_{2}\right) u^{2}}{64K%
}+-\frac{145Ku}{16}+\frac{61K_{1}^{2}-79K_{1}K_{2}+10K_{2}^{2}}{64K^{2}}u \\
&&+\frac{-17K_{11}+17K_{12}+13K_{22}}{32K}u-\frac{77\left(
K_{1}^{3}-2K_{1}^{2}K_{2}-K_{1}K_{2}^{2}+2K_{2}^{3}\right) }{64K^{3}} \\
&&+\frac{1}{96}\left( -365K_{1}+124K_{2}\right) -\frac{K_{2}\left(
32K_{11}+27K_{12}-95K_{22}\right) }{32K^{2}} \\
&&+\frac{K_{1}\left( 81K_{11}-79K_{12}-55K_{22}\right) }{64K^{2}}+\frac{%
-K_{111}+2K_{112}+K_{122}-3K_{222}}{4K};
\end{eqnarray*}

\begin{eqnarray*}
&&A_{120}= \\
&&\frac{3\left( K_{1}+K_{2}\right) u^{4}}{128K}-\frac{33\left(
K_{1}^{2}-K_{2}^{2}\right) u^{3}}{128K^{2}}+\frac{3\left(
K_{11}-K_{22}\right) u^{3}}{16K} \\
&&+\frac{23}{8}\left( K_{1}-K_{2}\right) u^{2}-\frac{23\left(
K_{1}+K_{2}\right) \left( K_{11}-K_{12}+K_{22}\right) u^{2}}{128K^{2}} \\
&&+\frac{5\left( K_{111}+K_{222}\right) u^{2}}{32K}-\frac{3\left(
33K_{1}^{2}-92K_{1}K_{2}+33K_{2}^{2}\right) u}{64K} \\
&&-\text{$\frac{77(K_{2}^{2}K_{11}+K_{1}^{2}K_{12})}{128K^{3}}u+\frac{%
5929K_{1}^{2}K_{2}^{2}K_{11}K_{12}u}{16384K^{6}}$}+\frac{31}{64}\left(
K_{11}-4K_{12}+K_{22}\right) u \\
&&-\frac{5\left( K_{11}-K_{22}\right) \left( K_{11}-K_{12}+K_{22}\right) u}{%
16K^{2}}+\frac{77\left( K_{1}^{2}-K_{2}^{2}\right) K_{22}u}{128K^{3}} \\
&&-\frac{K_{1}\left( 44K_{111}-15K_{112}+15K_{122}+6K_{222}\right) u}{%
128K^{2}} \\
&&+\frac{K_{2}\left( 6K_{111}+15K_{112}-15K_{122}+44K_{222}\right) }{128K^{2}%
}u+\frac{4K_{1111}+K_{1112}-K_{1222}-4K_{2222}}{32K}u \\
&&-\frac{11}{3}K\left( K_{1}+K_{2}\right) +\frac{17\left(
22K_{1}^{3}-53K_{1}^{2}K_{2}+53K_{1}K_{2}^{2}-22K_{2}^{3}\right) }{384K^{2}}
\\
&&+\frac{K_{1}\left( -1814K_{11}+402K_{12}-2011K_{22}\right) }{1920K}+\frac{%
K_{2}\left( 2011K_{11}-402K_{12}+1814K_{22}\right) }{1920K} \\
&&+\frac{77\left( K_{1}^{2}-K_{2}^{2}\right) \left( K_{112}-K_{122}\right) }{%
128K^{3}}-\frac{5\left( K_{11}-K_{22}\right) \left( K_{112}-K_{122}\right) }{%
16K^{2}} \\
&&-\frac{33}{40}\left( K_{111}-2K_{112}+2K_{122}-K_{222}\right) +\frac{%
K_{2}\left( 3K_{1112}+19K_{1122}-22K_{1222}\right) }{64K^{2}} \\
&&+\frac{K_{1}\left( -22K_{1112}+19K_{1122}+3K_{1222}\right) }{64K^{2}}+%
\frac{K_{11112}-K_{11122}-K_{11222}+K_{12222}}{8K},
\end{eqnarray*}

\begin{eqnarray*}
&&B_{121}= \\
&&-\frac{195u^{3}}{32}+\frac{9\left( -5K_{1}+9K_{2}\right) u^{2}}{16K}+\frac{%
169Ku}{16}+\frac{21\left( 26K_{1}^{2}-31K_{1}K_{2}+9K_{2}^{2}\right) u}{%
64K^{2}} \\
&&-\frac{3\left( 65K_{11}-75K_{12}+31K_{22}\right) u}{32K}-\frac{35K_{2}}{16}%
+\frac{49\left( 2K_{1}^{3}-3K_{1}^{2}K_{2}+3K_{1}K_{2}^{2}-K_{2}^{3}\right) 
}{32K^{3}} \\
&&+\frac{K_{1}\left( -58K_{11}+58K_{12}-37K_{22}\right) }{16K^{2}}+\frac{%
K_{2}\left( 29K_{11}-50K_{12}+29K_{22}\right) }{16K^{2}} \\
&&+\frac{6K_{111}-9K_{112}+9K_{122}-3K_{222}}{8K};
\end{eqnarray*}

\begin{eqnarray*}
&&B_{122}= \\
&&\frac{195u^{3}}{32}+\frac{9\left( 9K_{1}-5K_{2}\right) u^{2}}{16K}+\frac{%
169Ku}{16}-\frac{21\left( 9K_{1}^{2}-31K_{1}K_{2}+26K_{2}^{2}\right) u}{%
64K^{2}} \\
&&+\frac{3\left( 31K_{11}-75K_{12}+65K_{22}\right) u}{32K}+\frac{35K_{1}}{16}%
-\frac{49\left( K_{1}^{3}-3K_{1}^{2}K_{2}+3K_{1}K_{2}^{2}-2K_{2}^{3}\right) 
}{32K^{3}} \\
&&-\frac{K_{2}\left( 37K_{11}-58K_{12}+58K_{22}\right) }{16K^{2}}+\frac{%
K_{1}\left( 29K_{11}-50K_{12}+29K_{22}\right) }{16K^{2}} \\
&&-\frac{3\left( K_{111}-3K_{112}+3K_{122}-2K_{222}\right) }{8K};
\end{eqnarray*}

\begin{eqnarray*}
&&B_{120}= \\
&&\frac{15u^{5}}{64}+\frac{117\left( -K_{1}+K_{2}\right) u^{4}}{128K}-\frac{%
21\left( K_{1}^{2}-K_{1}K_{2}+K_{2}^{2}\right) u^{3}}{64K^{2}} \\
&&+\frac{65}{8}\left( K_{1}+K_{2}\right) u^{2}+\frac{273\left(
K_{1}-K_{2}\right) \left( K_{11}-K_{12}+K_{22}\right) u^{2}}{128K^{2}} \\
&&-\frac{3\left( 26K_{111}-47K_{112}+47K_{122}-26K_{222}\right) u^{2}}{64K}+%
\frac{33K^{2}u}{8} \\
&&+\frac{221\left( K_{1}^{2}-K_{2}^{2}\right) u}{64K}+\frac{351}{64}\left(
K_{11}-K_{22}\right) u-\frac{23\left( K_{11}-K_{12}+K_{22}\right) ^{2}u}{%
64K^{2}} \\
&&-\frac{7^{2}\left( K_{1}^{2}\left( K_{12}-K_{22}\right) -K_{2}^{2}\left(
K_{11}+K_{22}\right) +K_{1}K_{2}\left( K_{11}-K_{12}+K_{22}\right) \right) }{%
2^{6}K^{3}}u \\
&&-\frac{7^{4}K_{1}^{2}K_{2}^{2}K_{11}K_{12}}{2^{12}K^{6}}u-\frac{%
7K_{1}\left( 8K_{111}-51K_{112}+51K_{122}-4K_{222}\right) }{128K^{2}}u \\
&&-\frac{7K_{2}\left( -4K_{111}+51K_{112}-51K_{122}+8K_{222}\right) }{%
128K^{2}}u \\
&&+\frac{4K_{1111}-47K_{1112}+90K_{1122}-47K_{1222}+4K_{2222}}{32K}u \\
&&+\frac{119\left(
2K_{1}^{3}-3K_{1}^{2}K_{2}-3K_{1}K_{2}^{2}+2K_{2}^{3}\right) }{192K^{2}}+%
\frac{191K_{1}\left( -2K_{11}+2K_{12}+K_{22}\right) }{960K} \\
&&+\frac{191K_{2}\left( K_{11}+2K_{12}-2K_{22}\right) }{960K}+\frac{49\left(
K_{1}^{2}-K_{1}K_{2}+K_{2}^{2}\right) \left( K_{112}-K_{122}\right) }{64K^{3}%
} \\
&&-\frac{23\left( K_{11}-K_{12}+K_{22}\right) \left( K_{112}-K_{122}\right) 
}{64K^{2}}-\frac{33}{80}\left( 2K_{111}-3K_{112}-3K_{122}+2K_{222}\right) \\
&&-\frac{7\left( K_{1}\left( 2K_{1112}-3K_{1122}+K_{1222}\right) \right. }{%
32K^{2}}+\frac{K_{2}\left( K_{1112}-3K_{1122}+2K_{1222}\right) }{32K^{2}} \\
&&+\frac{K_{11112}-2K_{11122}+2K_{11222}-K_{12222}}{8K}.
\end{eqnarray*}

\endgroup

\subsection{The Polynomials $\emph{V}_{i}$ and $\emph{V}_{ij}$}

Here we give the expressions of the polynomials $V_{0},\;V_{1},\;V_{2},%
\;V_{11},\;V_{22},\;$and$\;V_{12}\;$(see Section 8) in terms of $%
A_{ij},\;A_{ijk},\;B_{ijk},\;$the curvature $K\;$and its first covariant
derivatives $K_{1}\;$and $K_{2}.$

\begingroup\scalefont{0.8}

\begin{eqnarray*}
&&V_{0}= \\
&&-12\left( K_{1}-2K_{2}\right) \left( 3K_{1}-K_{2}\right) A_{21}B_{122} \\
&&+12\left( 2K_{1}^{2}-7K_{1}K_{2}+3K_{2}^{2}\right) \left(
A_{11}B_{122}+A_{22}B_{121}-A_{12}B_{121}\right) \\
&&+39u\left( 11K_{1}^{2}-26K_{1}K_{2}+11K_{2}^{2}\right) \left(
A_{11}A_{22}-A_{12}A_{21}\right) \\
&&+16K\left( -2K_{1}+K_{2}\right) \left( A_{111}B_{122}-A_{112}B_{121}\right)
\\
&&+16K\left( K_{1}-2K_{2}\right) \left( A_{222}B_{121}-A_{221}B_{122}\right)
\\
&&+208Ku\left( K_{1}-2K_{2}\right) \left( A_{21}A_{222}-A_{22}A_{221}\right)
\\
&&+208Ku\left( 2K_{1}-K_{2}\right) \left( A_{11}A_{112}-A_{12}A_{111}\right)
\\
&&+52Ku\left( K_{1}+K_{2}\right) \left(
A_{11}A_{222}-A_{12}A_{221}-A_{21}A_{112}+A_{22}A_{111}\right) \\
&&+208K^{2}u\left( A_{111}A_{222}-A_{112}A_{221}\right) ;
\end{eqnarray*}

\begin{eqnarray*}
&&V_{1}= \\
&&12\left( 3K_{1}^{2}-7K_{1}K_{2}+2K_{2}^{2}\right) \left(
A_{20}B_{122}-A_{22}B_{120}\right) \\
&&+12\left( -2K_{1}^{2}+7K_{1}K_{2}-3K_{2}^{2}\right) \left(
A_{10}B_{122}-A_{12}B_{120}\right) \\
&&-39u\left( 11K_{1}^{2}-26K_{1}K_{2}+11K_{2}^{2}\right) \left(
A_{10}A_{22}-A_{12}A_{20}\right) \\
&&+16K\left( K_{1}-2K_{2}\right) \left( A_{220}B_{122}-A_{222}B_{120}\right)
\\
&&+16K\left( 2K_{1}-K_{2}\right) \left( A_{110}B_{122}-A_{112}B_{120}\right)
\\
&&-208Ku\left( K_{1}-2K_{2}\right) \left( A_{20}A_{222}-A_{22}A_{220}\right)
\\
&&-52Ku\left( K_{1}+K_{2}\right) \left(
A_{10}A_{222}-A_{12}A_{220}-A_{20}A_{112}+A_{22}A_{110}\right) \\
&&-208Ku\left( 2K_{1}-K_{2}\right) \left( A_{10}A_{112}-A_{12}A_{110}\right)
\\
&&-208K^{2}u\left( A_{110}A_{222}-A_{112}A_{220}\right) ;
\end{eqnarray*}

\begin{eqnarray*}
&&V_{2}= \\
&&-12\left( 3K_{1}^{2}-7K_{1}K_{2}+2K_{2}^{2}\right) \left(
A_{20}B_{121}-A_{21}B_{120}\right) \\
&&-12\left( 2K_{1}^{2}-7K_{1}K_{2}+3K_{2}^{2}\right) \left(
A_{11}B_{120}-A_{10}B_{121}\right) \\
&&+39u\left( 11K_{1}^{2}-26K_{1}K_{2}+11K_{2}^{2}\right) \left(
A_{10}A_{21}-A_{11}A_{20}\right) \\
&&-16K\left( K_{1}-2K_{2}\right) \left( A_{220}B_{121}-A_{221}B_{120}\right)
\\
&&-16K\left( 2K_{1}-K_{2}\right) \left( A_{110}B_{121}-A_{111}B_{120}\right)
\\
&&+208Ku\left( K_{1}-2K_{2}\right) \left( A_{20}A_{221}-A_{21}A_{220}\right)
\\
&&+52Ku\left( K_{1}+K_{2}\right) \left(
A_{10}A_{221}-A_{11}A_{220}-A_{20}A_{111}+A_{21}A_{110}\right) \\
&&+208Ku\left( 2K_{1}-K_{2}\right) \left( A_{10}A_{111}-A_{11}A_{110}\right)
\\
&&+208K^{2}u\left( A_{110}A_{221}-A_{111}A_{220}\right) ;
\end{eqnarray*}

\begin{eqnarray*}
&&V_{11}= \\
&&-\left( 52K_{1}^{2}-124K_{1}K_{2}+64K_{2}^{2}\right) \left(
A_{10}A_{21}B_{122}-A_{11}A_{20}B_{122}-A_{10}A_{22}B_{121}\right) \\
&&-\left( 52K_{1}^{2}-124K_{1}K_{2}+64K_{2}^{2}\right)
(A_{12}A_{20}B_{121}+A_{11}A_{22}B_{120}-A_{12}A_{21}B_{120}) \\
&&+\frac{16}{3}K\left( K_{1}-8K_{2}\right) \left(
A_{10}A_{221}B_{122}-A_{11}A_{220}B_{122}-A_{10}A_{222}B_{121}\right) \\
&&+\frac{16}{3}K\left( K_{1}-8K_{2}\right)
(A_{12}A_{220}B_{121}+A_{11}A_{222}B_{120}-A_{12}A_{221}B_{120}) \\
&&+\frac{32}{3}K\left( K_{1}+K_{2}\right) \left(
A_{20}A_{111}B_{122}-A_{21}A_{110}B_{122}-A_{20}A_{112}B_{121}\right) \\
&&+\frac{32}{3}K\left( K_{1}+K_{2}\right)
(A_{22}A_{110}B_{121}+A_{21}A_{112}B_{120}-A_{22}A_{111}B_{120}) \\
&&-\frac{80}{3}K\left( 2K_{1}-K_{2}\right) \left(
A_{10}A_{111}B_{122}-A_{11}A_{110}B_{122}-A_{10}A_{112}B_{121}\right) \\
&&-\frac{80}{3}K\left( 2K_{1}-K_{2}\right)
(+A_{12}A_{110}B_{121}+A_{11}A_{112}B_{120}-A_{12}A_{111}B_{120}) \\
&&+\frac{128}{3}K\left( K_{1}-2K_{2}\right) \left(
A_{20}A_{222}B_{121}-A_{22}A_{220}B_{121}-A_{21}A_{222}B_{120}\right) \\
&&-\frac{128}{3}K\left( K_{1}-2K_{2}\right) \left(
A_{20}A_{221}B_{122}-A_{21}A_{220}B_{122}-A_{22}A_{221}B_{120}\right) \\
&&-208Ku\left( K_{1}-2K_{2}\right) \left(
A_{10}A_{21}A_{222}-A_{11}A_{20}A_{222}-A_{10}A_{22}A_{221}\right) \\
&&-208Ku\left( K_{1}-2K_{2}\right)
(A_{12}A_{20}A_{221}+A_{11}A_{22}A_{220}-A_{12}A_{21}A_{220}) \\
&&+52Ku\left( K_{1}+K_{2}\right) \left(
A_{10}A_{21}A_{112}-A_{11}A_{20}A_{112}-A_{10}A_{22}A_{111}\right) \\
&&+52Ku\left( K_{1}+K_{2}\right)
(+A_{12}A_{20}A_{111}+A_{11}A_{22}A_{110}-A_{12}A_{21}A_{110}) \\
&&-\frac{128K^{2}}{3}\left(
A_{110}A_{221}B_{122}-A_{111}A_{220}B_{122}-A_{110}A_{222}B_{121}\right) \\
&&-\frac{128K^{2}}{3}%
(A_{112}A_{220}B_{121}+A_{111}A_{222}B_{120}-A_{112}A_{221}B_{120}) \\
&&-208K^{2}u(-A_{12}A_{111}A_{220}+A_{11}A_{112}A_{220}+A_{12}A_{110}A_{221})
\\
&&+208K^{2}u(A_{10}A_{112}A_{221}+A_{11}A_{110}A_{222}-A_{10}A_{111}A_{222});
\end{eqnarray*}

\begin{eqnarray*}
&&V_{22}= \\
&&-\left( 64K_{1}^{2}-124K_{1}K_{2}+52K_{2}^{2}\right) \left(
A_{10}A_{21}B_{122}-A_{11}A_{20}B_{122}-A_{10}A_{22}B_{121}\right) \\
&&-\left( 64K_{1}^{2}-124K_{1}K_{2}+52K_{2}^{2}\right)
(A_{12}A_{20}B_{121}+A_{11}A_{22}B_{120}-A_{12}A_{21}B_{120}) \\
&&-\frac{80}{3}K\left( K_{1}-2K_{2}\right) \left(
A_{20}A_{221}B_{122}-A_{21}A_{220}B_{122}-A_{20}A_{222}B_{121}\right) \\
&&-\frac{80}{3}K\left( K_{1}-2K_{2}\right)
(A_{22}A_{220}B_{121}+A_{21}A_{222}B_{120}-A_{22}A_{221}B_{120}) \\
&&-\frac{80}{3}K\left( K_{1}-2K_{2}\right)
(A_{12}A_{220}B_{121}+A_{11}A_{222}B_{120}-A_{12}A_{221}B_{120}) \\
&&-\frac{32}{3}K\left( K_{1}+K_{2}\right) \left(
A_{10}A_{221}B_{122}-A_{11}A_{220}B_{122}-A_{10}A_{222}B_{121}\right) \\
&&+\frac{16}{3}K\left( 8K_{1}-K_{2}\right) \left(
A_{20}A_{111}B_{122}-A_{21}A_{110}B_{122}-A_{20}A_{112}B_{121}\right) \\
&&+\frac{16}{3}K\left( 8K_{1}-K_{2}\right)
(+A_{22}A_{110}B_{121}+A_{21}A_{112}B_{120}-A_{22}A_{111}B_{120}) \\
&&-\frac{128}{3}K\left( 2K_{1}-K_{2}\right) \left(
A_{10}A_{111}B_{122}-A_{11}A_{110}B_{122}-A_{10}A_{112}B_{121}\right) \\
&&-\frac{128}{3}K\left( 2K_{1}-K_{2}\right)
(A_{12}A_{110}B_{121}+A_{11}A_{112}B_{120}-A_{12}A_{111}B_{120}) \\
&&+52Ku\left( K_{1}+K_{2}\right) \left(
A_{10}A_{21}A_{222}-A_{11}A_{20}A_{222}-A_{10}A_{22}A_{221}\right) \\
&&+52Ku\left( K_{1}+K_{2}\right)
(A_{12}A_{20}A_{221}+A_{11}A_{22}A_{220}-A_{12}A_{21}A_{220}) \\
&&+208Ku\left( 2K_{1}-K_{2}\right) \left(
A_{10}A_{21}A_{112}-A_{11}A_{20}A_{112}-A_{10}A_{22}A_{111}\right) \\
&&+208Ku\left( 2K_{1}-K_{2}\right)
(A_{12}A_{20}A_{111}+A_{11}A_{22}A_{110}-A_{12}A_{21}A_{110}) \\
&&-\frac{128K^{2}}{3}\left(
A_{110}A_{221}B_{122}-A_{111}A_{220}B_{122}-A_{110}A_{222}B_{121}\right) \\
&&-\frac{128K^{2}}{3}%
(A_{112}A_{220}B_{121}+A_{111}A_{222}B_{120}-A_{112}A_{221}B_{120}) \\
&&+208K^{2}u(A_{22}A_{111}A_{220}-A_{21}A_{112}A_{220}-A_{22}A_{110}A_{221})
\\
&&+208K^{2}u(+A_{20}A_{112}A_{221}+A_{21}A_{110}A_{222}-A_{20}A_{111}A_{222});
\end{eqnarray*}

\begin{eqnarray*}
&&V_{12}= \\
&&-\left( 44K_{1}^{2}-104K_{1}K_{2}+44K_{2}^{2}\right) \left(
A_{10}A_{21}B_{122}-A_{11}A_{20}B_{122}-A_{10}A_{22}B_{121}\right) \\
&&-\left( 44K_{1}^{2}-104K_{1}K_{2}+44K_{2}^{2}\right)
(A_{12}A_{20}B_{121}+A_{11}A_{22}B_{120}-A_{12}A_{21}B_{120}) \\
&&-\frac{64}{3}K\left( K_{1}-2K_{2}\right) \left(
A_{20}A_{221}B_{122}-A_{21}A_{220}B_{122}-A_{20}A_{222}B_{121}\right) \\
&&-\frac{64}{3}K\left( K_{1}-2K_{2}\right)
(A_{22}A_{220}B_{121}+A_{21}A_{222}B_{120}-A_{22}A_{221}B_{120}) \\
&&-\frac{64}{3}K\left( 2K_{1}-K_{2}\right) \left(
A_{10}A_{111}B_{122}-A_{11}A_{110}B_{122}-A_{10}A_{112}B_{121}\right) \\
&&-\frac{64}{3}K\left( 2K_{1}-K_{2}\right)
(A_{12}A_{110}B_{121}+A_{11}A_{112}B_{120}-A_{12}A_{111}B_{120}) \\
&&-\frac{16}{3}K\left( K_{1}+K_{2}\right)
(A_{10}A_{221}B_{122}-A_{11}A_{220}B_{122}-A_{20}A_{111}B_{122}) \\
&&-\frac{16}{3}K\left( K_{1}+K_{2}\right)
(A_{21}A_{110}B_{122}-A_{10}A_{222}B_{121}+A_{12}A_{220}B_{121}) \\
&&-\frac{16}{3}K\left( K_{1}+K_{2}\right)
(A_{20}A_{112}B_{121}-A_{22}A_{110}B_{121}+A_{11}A_{222}B_{120}) \\
&&+\frac{16}{3}K\left( K_{1}+K_{2}\right)
(A_{12}A_{221}B_{120}+A_{21}A_{112}B_{120}-A_{22}A_{111}B_{120}) \\
&&+\frac{64}{3}%
K^{2}(A_{112}A_{221}B_{120}-A_{111}A_{222}B_{120}-A_{112}A_{220}B_{121}) \\
&&+\frac{64}{3}%
K^{2}(+A_{110}A_{222}B_{121}+A_{111}A_{220}B_{122}-A_{110}A_{221}B_{122}).
\end{eqnarray*}

\endgroup

\subsection{ Total Covariant Derivatives}

\begingroup\scalefont{0.8}

\begin{eqnarray*}
&&\delta _{1}= \\
&&p_{1}\frac{\partial }{\partial u}+p_{11}\frac{\partial }{\partial p_{1}}%
+\left( -\frac{Ku}{2}+p_{12}\right) \frac{\partial }{\partial p_{2}}+p_{111}%
\frac{\partial }{\partial p_{11}} \\
&&+\left( p_{112}-\frac{uK_{1}}{6}-\frac{5Kp_{1}}{6}\right) \frac{\partial }{%
\partial p_{12}}+\left( p_{122}-\frac{uK_{2}}{3}-\frac{5Kp_{2}}{3}\right) 
\frac{\partial }{\partial p_{22}} \\
&&+p_{1111}\frac{\partial }{\partial p_{111}}+\left( p_{1112}-\frac{uK_{11}}{%
12}-\frac{K_{1}p_{1}}{2}-\frac{7Kp_{11}}{6}\right) \frac{\partial }{\partial
p_{112}} \\
&&+\left( p_{1122}-\frac{uK_{12}}{6}-\frac{K_{2}p_{1}}{2}-\frac{K_{1}p_{2}}{2%
}-\frac{7Kp_{12}}{3}\right) \frac{\partial }{\partial p_{122}} \\
&&+\left( p_{1222}-\frac{uK_{22}}{4}-\frac{3K_{2}p_{2}}{2}-\frac{7Kp_{22}}{2}%
\right) \frac{\partial }{\partial p_{222}};
\end{eqnarray*}

\begin{eqnarray*}
&&\delta _{2}= \\
&&p_{2}\frac{\partial }{\partial u}+\left( \frac{Ku}{2}+p_{12}\right) \frac{%
\partial }{\partial p_{1}}+p_{22}\frac{\partial }{\partial p_{2}}+\left(
p_{112}+\frac{uK_{1}}{3}+\frac{5Kp_{1}}{3}\right) \frac{\partial }{\partial
p_{11}} \\
&&+\left( p_{122}+\frac{uK_{2}}{6}+\frac{5Kp_{2}}{6}\right) \frac{\partial }{%
\partial p_{12}}+p_{222}\frac{\partial }{\partial p_{22}}+ \\
&&\left( p_{1112}+\frac{uK_{11}}{4}+\frac{3K_{1}p_{1}}{2}+\frac{7Kp_{11}}{2}%
\right) \frac{\partial }{\partial p_{111}} \\
&&+\left( p_{1122}+\frac{uK_{12}}{6}+\frac{K_{2}p_{1}}{2}+\frac{K_{1}p_{2}}{2%
}+\frac{7Kp_{12}}{3}\right) \ \frac{\partial }{\partial p_{112}} \\
&&+\left( p_{1222}+\frac{uK_{22}}{12}+\frac{K_{2}p_{2}}{2}+\frac{7Kp_{22}}{6}%
\right) \frac{\partial }{\partial p_{122}}+p_{2222}\frac{\partial }{\partial
p_{222}};
\end{eqnarray*}

\begin{eqnarray*}
&&\delta _{1}^{K}= \\
&&K_{1}\frac{\partial }{\partial K}+K_{11}\frac{\partial }{\partial K_{1}}%
+\left( -K^{2}+K_{12}\right) \frac{\partial }{\partial K_{2}}+K_{111}\frac{%
\partial }{\partial K_{11}} \\
&&+\left( K_{112}-\frac{5KK_{1}}{3}\right) \frac{\partial }{\partial K_{12}}%
+\left( K_{122}-\frac{10KK_{2}}{3}\right) \frac{\partial }{\partial K_{22}}
\\
&&+K_{1111}\frac{\partial }{\partial K_{111}}+\left( -\frac{5K_{1}^{2}}{6}-%
\frac{11KK_{11}}{6}+K_{1112}\right) \frac{\partial }{\partial K_{112}} \\
&&+\left( K_{1122}-\frac{5}{3}K_{1}K_{2}-\frac{11KK_{12}}{3}\right) \frac{%
\partial }{\partial K_{122}}+\left( K_{1222}-\frac{5K_{2}^{2}}{2}-\frac{%
11KK_{22}}{2}\right) \frac{\partial }{\partial K_{222}} \\
&&+K_{11111}\frac{\partial }{\partial K_{1111}}+\left( K_{11112}-\frac{%
21K_{1}K_{11}}{10}-\frac{21KK_{111}}{10}\right) \frac{\partial }{\partial
K_{1112}} \\
&&+\left( K_{11122}-\frac{7K_{2}K_{11}}{5}-\frac{14K_{1}K_{12}}{5}-\frac{%
21KK_{112}}{5}\right) \frac{\partial }{\partial K_{1122}} \\
&&+\left( K_{11222}-\frac{21K_{2}K_{12}}{5}-\frac{21K_{1}K_{22}}{10}-\frac{%
63KK_{122}}{10}\right) \frac{\partial }{\partial K_{1222}} \\
&&+\left( K_{12222}-\frac{42K_{2}K_{22}}{5}-\frac{42KK_{222}}{5}\right) 
\frac{\partial }{\partial K_{2222}}+K_{111111}\frac{\partial }{\partial
K_{11111}} \\
&&+\left( K_{111112}-\frac{7K_{11}^{2}}{5}-\frac{14K_{1}K_{111}}{5}-\frac{%
12KK_{1111}}{5}\right) \frac{\partial }{\partial K_{11112}} \\
&&+\left( K_{111122}-\frac{14K_{11}K_{12}}{5}-\frac{7K_{2}K_{111}}{5}-\frac{%
21K_{1}K_{112}}{5}-\frac{24KK_{1112}}{5}\right) \frac{\partial }{\partial
K_{11122}} \\
&&+\left( K_{111222}-\frac{14K_{12}^{2}}{5}-\frac{7K_{11}K_{22}}{5}-\frac{%
21K_{2}K_{112}}{5}-\frac{21K_{1}K_{122}}{5}-\frac{36KK_{1122}}{5}\right) 
\frac{\partial }{\partial K_{11222}} \\
&&+\left( K_{112222}-\frac{28K_{12}K_{22}}{5}-\frac{42K_{2}K_{122}}{5}-\frac{%
14K_{1}K_{222}}{5}-\frac{48KK_{1222}}{5}\right) \frac{\partial }{\partial
K_{12222}} \\
&&+\left( K_{122222}-7K_{22}^{2}-14K_{2}K_{222}-12KK_{2222}\right) \frac{%
\partial }{\partial K_{22222}};
\end{eqnarray*}

\begin{eqnarray*}
&&\delta _{2}^{K}= \\
&&K_{2}\frac{\partial }{\partial K}+\left( K_{12}+K^{2}\right) \frac{%
\partial }{\partial K_{1}}+K_{22}\frac{\partial }{\partial K_{2}}+\left(
K_{112}+\frac{10KK_{1}}{3}\right) \frac{\partial }{\partial K_{11}} \\
&&+\left( K_{122}+\frac{5KK_{2}}{3}\right) \frac{\partial }{\partial K_{12}}%
+\left( K_{1112}+\frac{5K_{1}^{2}}{2}+\frac{11KK_{11}}{2}\right) \frac{%
\partial }{\partial K_{111}} \\
&&+K_{222}\frac{\partial }{\partial K_{22}}+\left( K_{1122}+\frac{5K_{1}K_{2}%
}{3}+\frac{11KK_{12}}{3}\right) \frac{\partial }{\partial K_{112}}+K_{2222}%
\frac{\partial }{\partial K_{222}} \\
&&+\left( K_{1222}+\frac{5K_{2}^{2}}{6}+\frac{11KK_{22}}{6}\right) \frac{%
\partial }{\partial K_{122}}+K_{22222}\frac{\partial }{\partial K_{2222}} \\
&&+\left( K_{12222}+\frac{21K_{2}K_{22}}{10}+\frac{21KK_{222}}{10}\right) 
\frac{\partial }{\partial K_{1222}} \\
&&+\left( K_{11222}+\frac{14K_{2}K_{12}}{5}+\frac{7K_{1}K_{22}}{5}+\frac{%
21KK_{122}}{5}\right) \frac{\partial }{\partial K_{1122}} \\
&&+\left( K_{11122}+\frac{21K_{2}K_{11}}{10}+\frac{21K_{1}K_{12}}{5}+\frac{%
63KK_{112}}{10}\right) \frac{\partial }{\partial K_{1112}} \\
&&+\left( K_{11112}+\frac{42K_{1}K_{11}}{5}+\frac{42KK_{111}}{5}\right) 
\frac{\partial }{\partial K_{1111}} \\
&&+K_{222222}\frac{\partial }{\partial K_{22222}}+\left( K_{122222}+\frac{%
7K_{22}^{2}}{5}+\frac{14K_{2}K_{222}}{5}+\frac{12KK_{2222}}{5}\right) \frac{%
\partial }{\partial K_{12222}} \\
&&+\left( K_{112222}+\frac{14K_{12}K_{22}}{5}+\frac{21K_{2}K_{122}}{5}+\frac{%
7K_{1}K_{222}}{5}+\frac{24KK_{1222}}{5}\right) \frac{\partial }{\partial
K_{11222}} \\
&&+\left( K_{111222}+\frac{14K_{12}^{2}}{5}+\frac{7K_{11}K_{22}}{5}+\frac{%
21K_{2}K_{112}}{5}+\frac{21K_{1}K_{122}}{5}+\frac{36KK_{1122}}{5}\right) 
\frac{\partial }{\partial K_{11122}} \\
&&+\left( K_{111122}+\frac{28K_{11}K_{12}}{5}+\frac{14K_{2}K_{111}}{5}+\frac{%
42K_{1}K_{112}}{5}+\frac{48KK_{1112}}{5}\right) \frac{\partial }{\partial
K_{11112}} \\
&&+\left( K_{111112}+7K_{11}^{2}+14K_{1}K_{111}+12KK_{1111}\right) \frac{%
\partial }{\partial K_{11111}}.
\end{eqnarray*}

\endgroup

{\emph{Authors' addresses:} }

{\ Deparment of Mathematical Sciences, New Jersey Institute of Technology,
University Heights, Newark, NJ 07102, USA; vlgold@oak.njit.edu }

{\ Department of Mathematics, The University of Tromso, N9037, Tromso,
Norway; lychagin@math.uit.no }

\end{document}